\documentclass[reqno]{amsart}

\usepackage{amssymb,mathtools,mathrsfs,bm,enumitem,booktabs}

\usepackage[hidelinks]{hyperref}
\hypersetup{
  pdftitle={Exact-Order Mellin Transfer and Prime-Support Filtrations},
  pdfauthor={Dongkyu Lim; Abdelmejid Bayad},
  pdfkeywords={Mellin transforms, Mobius inversion, prime-support filtrations, spectral zeta functions, mod-Poisson convergence, Galois entanglement, exact-order torsion}
}
\usepackage[nameinlink,capitalise,noabbrev]{cleveref}

\numberwithin{equation}{section}
\allowdisplaybreaks[4]
\theoremstyle{plain}
\newtheorem{theorem}{Theorem}[section]
\newtheorem{proposition}[theorem]{Proposition}
\newtheorem{corollary}[theorem]{Corollary}
\newtheorem{lemma}[theorem]{Lemma}

\theoremstyle{definition}

\newtheorem{example}[theorem]{Example}
\newtheorem{problem}[theorem]{Problem}

\theoremstyle{remark}
\newtheorem{remark}[theorem]{Remark}

\newcommand{\R}{\mathbb R}
\newcommand{\C}{\mathbb C}
\newcommand{\Z}{\mathbb Z}
\newcommand{\E}{\mathbb E}
\newcommand{\Prob}{\mathbb P}
\newcommand{\ba}{\bm a}
\newcommand{\bb}{\bm b}

\newcommand{\BB}{\mathcal B}
\newcommand{\Ca}{\mathcal C}
\newcommand{\Fa}{\mathcal F}
\newcommand{\Ta}{\mathcal T}
\newcommand{\Ra}{\mathcal R}

\newcommand{\Xa}{X_{\ba}}
\newcommand{\dd}{\,\mathrm d}

\newcommand{\ord}{\operatorname{ord}}
\newcommand{\rad}{\operatorname{rad}}

\begin{document}

\title[Exact-order Mellin transfer and prime-support filtrations]
{Exact-Order Mellin Transfer and Prime-Support Filtrations}

\author[D. Lim]{Dongkyu Lim}
\address{Department of Mathematics Education, Gyeongkuk National University, Andong 36729, Republic of Korea}
\email{dklim@gknu.ac.kr}
\urladdr{\url{https://orcid.org/0000-0002-0928-8480}}

\author[A. Bayad]{Abdelmejid Bayad}
\address{Laboratoire de Math\'ematiques et Mod\'elisation d'Evry (LaMME), CNRS (UMR 8071), Universit\'e Paris-Saclay, Evry cedex, France}
\email{abdelmejid.bayad@univ-evry.fr}
\urladdr{\url{https://orcid.org/0000-0002-8003-335X}}
\begin{abstract}
We establish a unified exact-order Mellin transfer principle that converts the
multiplicity of a critical Mellin pole into arithmetic prime-support depth.
For an exact-order M\"obius transform in lattice dimension $\ell$, a base pole
of order $r$ at $s=\ell/2$ satisfies
\[
 \ord_{s=\ell/2}\widehat{\mathscr P_{N,\ell}H}=\omega(N)-r.
\]
The complete Jordan--Mellin jet is explicit, and its first surviving
coefficient selects the canonical weight $W(N)=\prod_{p\mid N}\log p$.
After the corresponding jet normalization, this weight governs a
kernel-independent moving-support law and a mod-Poisson regime.

With simultaneous splitting marks in finite Galois extensions, the same
weight yields a multivariate polynomial whose mixed coefficients are
reciprocal compositum degrees; it recovers pairwise intersection degrees and,
under a full-rank hypothesis, also governs the coupled saddle and
large-deviation geometry.  In spectral geometry, a two-variable family
attached to any positive self-adjoint elliptic operator has a transverse
normal crossing whose parameter jets generate every higher critical pole
order; exact-order dilation therefore realizes every positive prime-support
depth.  An explicit Bernoulli--Barnes/Barnes--Epstein model shows how a pole
collision shifts the finite critical layer from one-prime to two-prime
support.  Thus arithmetic support, Galois compositum data, and spectral pole
geometry arise from one exact-order Mellin mechanism.
\end{abstract}


\thanks{Corresponding author: Abdelmejid Bayad.}
\maketitle
\newpage
\tableofcontents
\newpage

\section{Introduction}

Mellin singularities usually encode analytic or geometric data, whereas
M\"obius inversion separates exact arithmetic order.  The main point of this
paper is that, at the critical Mellin point, these two operations obey one
rigid transfer mechanism: exact-order inversion contributes a zero whose
order is the number of distinct prime divisors, and that zero cancels the
critical pole multiplicity one order at a time.  The resulting stratification
is the \emph{prime-support filtration}.

For a Mellin kernel $H$ in lattice dimension $\ell$, put
\[
 (\mathscr P_{N,\ell}H)(t)
 =\sum_{d\mid N}\mu_{\rm M}(N/d)d^\ell H(d^2t).
\]
Mellin transformation gives
\[
 \widehat{\mathscr P_{N,\ell}H}(s)
 =J_{\ell-2s}(N)\widehat H(s).
\]
At $s_\ell=\ell/2$, the Jordan factor has a zero of exact order
$q=\omega(N)$ and first nonzero Taylor coefficient
$q!\prod_{p\mid N}\log p$.  Hence, if $\widehat H$ has a pole of order $r$
at $s_\ell$,
\begin{equation}\label{eq:intro-transfer-law-v59}
 \ord_{s=s_\ell}\widehat{\mathscr P_{N,\ell}H}=\omega(N)-r.
\end{equation}
The complete Laurent--Taylor convolution is equally important: its first
surviving coefficient selects the canonical weight
\[
 W(N):=\prod_{p\mid N}\log p,
\]
and the same weight reappears in the global, Galois, and spectral parts of
the paper.

The central theme is that local pole multiplicity, moving prime support,
Frobenius marking, and spectral pole geometry are different manifestations of
the same exact-order transfer law.  The following theorem records the four
components of this mechanism; the later sections establish them in a single
logical progression.

\begin{theorem}\label{thm:unified-main}
The following statements hold.
\begin{enumerate}[label=\textup{(\Roman*)},leftmargin=3em]
\item \textbf{Local exact-order transfer.}
Let $\widehat H$ be meromorphic near $s_\ell=\ell/2$ with a pole of order
$r\ge1$ and leading Laurent coefficient $c_{-r}\ne0$.  If
$q=\omega(N)$, then
\[
 \ord_{s=s_\ell}\widehat{\mathscr P_{N,\ell}H}=q-r,
\]
and, for $q\ge r$,
\[
 \left.\frac{\dd^{q-r}}{\dd s^{q-r}}
 \widehat{\mathscr P_{N,\ell}H}(s)\right|_{s=s_\ell}
 =\frac{(q-r)!(-2)^q}{q!}\Lambda_q(N)c_{-r},
 \qquad
 \Lambda_q(N)=q!\prod_{p\mid N}\log p.
\]
The full critical jet is the Laurent--Taylor convolution of the Mellin germ
with the generalized von Mangoldt jet of the Jordan factor.

\item \textbf{Moving-support universality after jet normalization.}
For
\[
 \mathcal F(\sigma,z)
 =\sum_{N\ge1}\frac{W(N)z^{\omega(N)}}{N^\sigma},
 \qquad W(N)=\prod_{p\mid N}\log p,
\]
there is an entire function $\mathcal C_0$ such that, locally uniformly in
$z$,
\[
 e^{-z/\varepsilon}\mathcal F(1+\varepsilon,z)
 =\mathcal C_0(z)+O(\varepsilon),
\]
and, uniformly for $q\varepsilon$ in compact sets,
\[
 q!\varepsilon^q
 \sum_{\omega(N)=q}\frac{W(N)}{N^{1+\varepsilon}}
 =\mathcal C_0(q\varepsilon)+O(\varepsilon).
\]
The induced support variable has a mod-Poisson limit with parameter
$1/\varepsilon$, together with local asymptotics and a Berry--Esseen bound of
order $O(\sqrt\varepsilon)$.  After the natural normalization, the same
$\mathcal C_0$ governs the first surviving jets of every kernel in part~\textup{(I)}.

\item \textbf{Frobenius and Galois refinement.}
A one-dimensional character twist replaces $\omega(N)$ by the number of
selected prime divisors.  More generally, let $L_1,\ldots,L_r/\mathbb Q$ be
finite Galois extensions, let $M_i$ count prime divisors that split completely
in $L_i$, and let $L_S$ be the compositum of the fields indexed by
$\emptyset\ne S\subseteq[r]$.  Then the marked boundary law is governed by
\[
 P_{\mathbf L}(\mathbf z)
 =1+\sum_{\emptyset\ne S\subseteq[r]}
   \frac1{[L_S:\mathbb Q]}\prod_{i\in S}(z_i-1).
\]
Every fixed mixed cumulant with support $S$ satisfies
\[
 \kappa(M_{i,\varepsilon}:i\in S)
 =\frac1{\varepsilon[L_S:\mathbb Q]}+O(1),
\]
so pairwise normalized covariance recovers
$[L_i\cap L_j:\mathbb Q]$.  Under the full-rank marking hypothesis, the same
polynomial determines the coupled saddle equations, the Gaussian covariance,
and the large-deviation rate function.  The prime-sum input is obtained from
Dedekind zeta functions of the composita and is unconditional.  Here
``entanglement'' refers to the arithmetic coupling imposed by the common
composita on the mixed coefficients and support lattice; it does not assert
primitive probabilistic dependence among the disjoint exact Frobenius-pattern
classes.

\item \textbf{Spectral hierarchy and geometric realization.}
Let $P$ be a positive self-adjoint classical elliptic operator of order
$\rho>0$ on a closed $n$-manifold and put $a=n/\rho$.  For $u$ sufficiently
near $0$ and $t>0$, the operator difference inside
\[
 \mathfrak S_P(u;t)
 :=\operatorname{Tr}'\!\left((I+tP)^{-(a-u)}-(tP)^{-(a-u)}\right)
\]
is trace class, and the scalar family has Mellin transform
\[
 \Gamma(s)\frac{\Gamma(a-u-s)}{\Gamma(a-u)}\zeta_P(s),
\]
whose polar germ at $(u,s)=(0,a)$ is a transverse normal crossing.  Its
$j$th parameter jet has a pole of exact order $j+2$, and exact-order dilation
therefore gives critical valuation
\[
 \omega(N)-(j+2).
\]
Thus one spectral datum realizes every positive prime-support depth; on flat
tori this is exact-order torsion summation.  In the explicit
Bernoulli--Barnes/Barnes--Epstein model, collision of the Epstein pole with a
normalized-remainder pole raises the critical multiplicity from one to two
exactly in dimensions $2,4,\ldots,2m+2$, shifting the finite nonzero layer
from prime powers to integers supported on two distinct primes.
\end{enumerate}
\end{theorem}

The four parts are established by the following results.  Part~\textup{(I)} is
Theorem~\ref{thm:AMR-local}.  Part~\textup{(II)} follows from
Theorems~\ref{thm:moving-support} and~\ref{thm:mod-Poisson-jet}, together
with Corollary~\ref{cor:moving-first-jets}.  Part~\textup{(III)} is established
first in one dimension by Theorems~\ref{thm:twisted-local-resonance} and
\ref{thm:Frobenius-moving-support}, and then in its multivariate form by
Theorems~\ref{thm:Galois-entanglement-support} and~\ref{thm:Galois-moving-profile}.
Part~\textup{(IV)} combines the explicit pole-collision transition of
Theorem~\ref{thm:dimension-arithmetic-v14}, the general spectral normal-crossing
theory in Theorem~\ref{thm:spectral-master-family} and
Corollary~\ref{cor:complete-spectral-hierarchy}, and the flat-torus realization
of Theorem~\ref{thm:flat-torus-hierarchy}.  Thus every component of the unified
statement is tied to a specific theorem below, while the general spectral result
serves as the structural culmination rather than an unrelated application.

\paragraph{Relation to existing work and precise novelty.}
The analytic inputs are classical: elliptic complex powers and spectral zeta
functions \cite{Seeley1967,Hormander1968}; meromorphic pseudodifferential
families \cite{Lesch1999,Grubb2005,PaychaScott2007,AzzaliEtAl2015}; Barnes,
Epstein, and Kronecker theory \cite{Barnes1904,BayadBeck2014,Epstein1903,Kronecker1857,Lang1987};
Dedekind zeta functions \cite{Neukirch1999}; and established results on
Frobenian prime sets, restricted prime factors, and mod-$\phi$/mod-Poisson
convergence
\cite{Serre2011,LoughranMatthiesen2024,Norton1979,Mangerel2016,FerayMeliotNikeghbali2016,KowalskiNikeghbali2010}.

The new point is the exact-order transfer mechanism linking these inputs.
Exact-order summation does more than isolate arithmetic order: its complete
critical Jordan jet records pole multiplicity and selects $W(N)$, from which
the normalized moving boundary law follows.  In the Galois setting, the marked
law has coefficients $[L_S:\mathbb Q]^{-1}$, recovers compositum and
intersection data, and uses the same polynomial for its saddle and
large-deviation geometry.  Spectrally, parameter jets of a transverse normal
crossing generate higher pole orders, and exact-order dilation converts them
into the complete prime-support hierarchy.  The Bernoulli--Barnes/Barnes--Epstein
model supplies an explicit collision in which the finite critical layer moves
from one-prime to two-prime support.

\paragraph{Organization.}
Section~\ref{sec:abstract-resonance} proves the local transfer theorem and its
fixed-support lift.  Section~\ref{sec:prime-support-statistics} develops the normalized moving-support and mod-Poisson laws.  Section~\ref{sec:twisted-resonance}
proves the character and Galois refinements.  Sections~\ref{sec:centered}--\ref{sec:higher-dimensional-resonance}
construct the Bernoulli--Barnes/Barnes--Epstein realization and the explicit
dimension--arithmetic transition.  Section~\ref{sec:spectral-hierarchy} then
proves the general spectral normal-crossing hierarchy, and
Section~\ref{sec:flat-tori} gives its flat-torus exact-order torsion
realization.  The final section records the structural conclusions and open
problems.

\paragraph{Notation guide.}
Three exact-order operators are used: $\mathscr P_{N,\ell}$ for the basic
lattice scaling, $\mathscr P_{N,\ell,\chi}$ for its character twist, and
$\mathscr P^{(\rho)}_{N,n}$ for general dilation order $\rho$.  The arithmetic
quantities are $\omega(N)$, $W(N)=\prod_{p\mid N}\log p$, the Jordan factor
$J_\nu(N)$, and its Taylor coefficients $\Lambda_k(N)$.  In the multivariate
Galois section, $P_{\mathbf L}$ is the compositum polynomial,
$\mathcal V_{\mathbf L}$ the active jump set, $\Lambda_{\mathbf L}$ its
integer span, $\Gamma_{\mathbf L}$ the generated semigroup, and
$\mathscr C_{\mathbf L}$ the associated cone.  These groups of notation are
kept separate throughout.

\section{Exact-order Mellin resonance}
\label{sec:abstract-resonance}

We establish here the local mechanism underlying all subsequent results.
Exact-order M\"obius inversion introduces a Jordan factor whose vanishing
order at the critical Mellin point is exactly the number of distinct prime
divisors of $N$.  When this factor is combined with a Mellin transform having
a pole of order $r$, the resulting critical order is $\omega(N)-r$.
Theorem~\ref{thm:AMR-local} gives not only this valuation formula but the
complete Laurent--Taylor jet, with generalized von Mangoldt functions as its
arithmetic coefficients.  Theorem~\ref{thm:AMR-local} is the
structural engine of the paper: beyond the Mellin factorization, exact-order
averaging creates a zero whose order is exactly the number of distinct prime
divisors of $N$, and this zero offsets the critical pole order of the base
Mellin transform.

\subsection{Local resonance and generalized von Mangoldt jets}

For $\ell\in\mathbb N$ and $N\ge1$ define
\begin{equation}\label{eq:abstract-P}
 (\mathscr P_{N,\ell}H)(t)
 :=\sum_{d\mid N}\mu_{\rm M}(N/d)d^\ell H(d^2t).
\end{equation}
Here $\mu_{\rm M}$ denotes the classical M\"obius function \cite{Mobius1832}.
Let
\[
 J_\nu(N):=\sum_{d\mid N}\mu_{\rm M}(N/d)d^\nu,
 \qquad
 \Lambda_k(N):=\sum_{d\mid N}\mu_{\rm M}(N/d)(\log d)^k.
\]
The function $J_\nu$ is the complex-exponent analogue of the classical Jordan totient factor; see \cite{Jordan1870}.
Thus $\Lambda_1$ is the classical von Mangoldt function introduced in \cite{vonMangoldt1894};
we write $\Lambda(N):=\Lambda_1(N)$.  The functions $\Lambda_k$ themselves are classical; see, for example,
\cite{Ivic1985}.  What is decisive here is that exact-order M\"obius
inversion forces them to occur as the complete critical Taylor jet of the
Jordan factor in \eqref{eq:abstract-Jordan-jet}.  We also put
\[
 \omega(N):=\#\{p:p\mid N\},\qquad
 \rad(N):=\prod_{p\mid N}p,
\]
with $\omega(1)=0$, and put $s_\ell:=\ell/2$ and use the convention that $\ord_{s=s_\ell}f$ is positive
for a zero and negative for a pole.

\begin{theorem}\label{thm:AMR-local}
Assume that $H$ has an ordinary Mellin transform on a nonempty vertical
strip and that this transform admits a meromorphic continuation to a
neighborhood of $s_\ell$ with
\begin{equation}\label{eq:abstract-H-Laurent}
 \widehat H(s_\ell+u)=\sum_{j=-r}^{\infty}c_j u^j,
 \qquad r\ge1,\quad c_{-r}\ne0.
\end{equation}
Then, first on every common strip of absolute convergence and afterwards by
meromorphic continuation,
\begin{equation}\label{eq:abstract-Jordan-factor}
 \widehat{\mathscr P_{N,\ell}H}(s)
 =J_{\ell-2s}(N)\widehat H(s).
\end{equation}
Moreover,
\begin{equation}\label{eq:abstract-Jordan-jet}
 J_{\ell-2(s_\ell+u)}(N)
 =\sum_{k=0}^{\infty}\frac{(-2)^k}{k!}\Lambda_k(N)u^k.
\end{equation}
If
\[
 \widehat{\mathscr P_{N,\ell}H}(s_\ell+u)
 =\sum_{m=-r}^{\infty}C_m(N)u^m,
\]
then the complete jet is the convolution
\begin{equation}\label{eq:abstract-full-convolution}
 C_m(N)=\sum_{k=0}^{m+r}\frac{(-2)^k}{k!}\Lambda_k(N)c_{m-k}.
\end{equation}
For $N\ge2$, let $q=\omega(N)$.  Then
\begin{equation}\label{eq:abstract-support}
 \Lambda_k(N)=0\quad(k<q),
 \qquad
 \Lambda_q(N)=q!\prod_{p\mid N}\log p,
\end{equation}
and consequently
\begin{equation}\label{eq:abstract-valuation}
 \ord_{s=s_\ell}\widehat{\mathscr P_{N,\ell}H}=q-r.
\end{equation}
If $q\ge r$, the first surviving derivative is
\begin{equation}\label{eq:abstract-first-jet}
 \left.\frac{\dd^{q-r}}{\dd s^{q-r}}
 \widehat{\mathscr P_{N,\ell}H}(s)\right|_{s=s_\ell}
 =\frac{(q-r)!(-2)^q}{q!}\Lambda_q(N)c_{-r}.
\end{equation}
\end{theorem}

\begin{proof}
Fix a vertical strip on which all Mellin integrals below converge absolutely.
For each divisor $d\mid N$, the change of variables $x=d^2t$ gives
\[
 \int_0^\infty t^{s-1}d^\ell H(d^2t)\dd t
 =d^{\ell-2s}\int_0^\infty x^{s-1}H(x)\dd x
 =d^{\ell-2s}\widehat H(s).
\]
Since the divisor sum is finite, it can be interchanged with the integral,
and summation over $d\mid N$ gives
\eqref{eq:abstract-Jordan-factor}.  Both sides are meromorphic near $s_\ell$,
so the identity continues there.

Write $s=s_\ell+u$.  Since $\ell-2s=-2u$,
\[
 J_{-2u}(N)
 =\sum_{d\mid N}\mu_{\rm M}(N/d)e^{-2u\log d}
 =\sum_{k\ge0}\frac{(-2u)^k}{k!}
   \sum_{d\mid N}\mu_{\rm M}(N/d)(\log d)^k,
\]
which is \eqref{eq:abstract-Jordan-jet}.  Multiplying this convergent Taylor
series by the Laurent expansion \eqref{eq:abstract-H-Laurent} and collecting
the coefficient of $u^m$ gives \eqref{eq:abstract-full-convolution}.

It remains to identify the first nonzero Jordan coefficient.  If
$N=\prod_{j=1}^q p_j^{a_j}$, then
\[
 J_\nu(N)=N^\nu\prod_{p\mid N}(1-p^{-\nu}).
\]
Hence, with $\nu=-2u$,
\[
 N^{-2u}=1+O(u),\qquad
 1-p^{2u}=-2u\log p+O(u^2),
\]
and therefore
\begin{equation}\label{eq:Jordan-leading-local}
 J_{-2u}(N)
 =(-2)^q\!\left(\prod_{p\mid N}\log p\right)u^q
  +O(u^{q+1}).
\end{equation}
Comparison with \eqref{eq:abstract-Jordan-jet} proves
$\Lambda_k(N)=0$ for $k<q$ and
$\Lambda_q(N)=q!\prod_{p\mid N}\log p$.
Finally,
\[
 \widehat H(s_\ell+u)=c_{-r}u^{-r}+O(u^{-r+1}),
\]
so \eqref{eq:Jordan-leading-local} yields
\[
 \widehat{\mathscr P_{N,\ell}H}(s_\ell+u)
 =(-2)^q c_{-r}\!\left(\prod_{p\mid N}\log p\right)
 u^{q-r}+O(u^{q-r+1}).
\]
This proves the valuation formula.  When $q\ge r$, multiplying the leading
coefficient by $(q-r)!$ and using
$\Lambda_q(N)=q!\prod_{p\mid N}\log p$ gives
\eqref{eq:abstract-first-jet}.
\end{proof}

\begin{remark}\label{rem:critical-support-rigidity}
At the critical layer $q=r$, the critical value depends on the prime
support $\rad(N)$ and not on prime exponents.  Higher jets involve the
full convolution \eqref{eq:abstract-full-convolution}; only the first
surviving coefficient is purely $\Lambda_q(N)$.
\end{remark}

\subsection{Marked global lift and fixed-support scaling}

The local resonance law persists after global summation.  We introduce a
marked Dirichlet series that records the prime-support depth and show that,
for fixed support size, the same exponent $q-r$ governing the local critical
order reappears in the global boundary profile.  Thus the local valuation is
not lost under averaging; it survives as a global scaling law.  For
$\nu=\ell-2s$ define
\begin{equation}\label{eq:abstract-Euler}
 \mathscr E_\ell(w,s;z)
 :=\sum_{N\ge1}\frac{z^{\omega(N)}J_\nu(N)}{N^w}
 =\prod_p\left(
 1+z\frac{(p^\nu-1)p^{-w}}{1-p^{\nu-w}}
 \right)
\end{equation}
in the absolute-convergence region
$\Re w>1$ and $\Re(w-\nu)>1$.  At $z=1$,
\begin{equation}\label{eq:abstract-zeta-quotient}
 \mathscr E_\ell(w,s;1)
 =\frac{\zeta(w+2s-\ell)}{\zeta(w)}.
\end{equation}

\begin{theorem}\label{thm:AMR-fixedq}
Under the hypotheses of Theorem~\ref{thm:AMR-local}, let $q\ge0$ and
$\lambda>0$ be fixed.  Then
\begin{equation}\label{eq:abstract-fixedq}
 \lim_{\varepsilon\downarrow0}\varepsilon^r
 \sum_{\substack{N\ge1\\\omega(N)=q}}
 \frac{\widehat{\mathscr P_{N,\ell}H}
 (s_\ell+\lambda\varepsilon)}{N^{1+\varepsilon}}
 =c_{-r}\lambda^{-r}
 \frac{[-\log(1+2\lambda)]^q}{q!}.
\end{equation}
In particular, as $\lambda\downarrow0$, the right side is
\begin{equation}\label{eq:abstract-local-global}
 \frac{(-2)^q c_{-r}}{q!}\lambda^{q-r}
 +O(\lambda^{q-r+1}).
\end{equation}
Thus the same exponent $q-r$ occurs in the local valuation and the global
fixed-support critical profile.
\end{theorem}

\begin{proof}
The local factor in \eqref{eq:abstract-Euler} is obtained from
$J_\nu(p^a)=p^{a\nu}-p^{(a-1)\nu}$ by summing over $a\ge1$.
At $z=1$ this factor equals
$(1-p^{-w})/(1-p^{\nu-w})$, proving
\eqref{eq:abstract-zeta-quotient}.

Put $w=1+\varepsilon$ and $s=s_\ell+\lambda\varepsilon$, and set
\[
 A_p(\varepsilon,\lambda)
 :=\frac{(p^{-2\lambda\varepsilon}-1)p^{-1-\varepsilon}}
 {1-p^{-1-(1+2\lambda)\varepsilon}}.
\]
Then the Euler product in \eqref{eq:abstract-Euler} is simply
$\prod_p(1+zA_p)$.  At $z=1$, \eqref{eq:abstract-zeta-quotient} gives
\[
 \prod_p(1+A_p)
 =\frac{\zeta(1+(1+2\lambda)\varepsilon)}{\zeta(1+\varepsilon)}.
\]
Consequently
\begin{equation}\label{eq:marked-factorization-detail}
 \mathscr E_\ell(1+\varepsilon,s_\ell+\lambda\varepsilon;z)
 =\left[
 \frac{\zeta(1+(1+2\lambda)\varepsilon)}
 {\zeta(1+\varepsilon)}\right]^z
 \prod_p\frac{1+zA_p}{(1+A_p)^z}.
\end{equation}
For $\lambda$ in a fixed compact subset of $(0,\infty)$,
$|1-e^{-x}|\ll\min\{x,1\}$ and the denominator is bounded below by
$1-p^{-1}$, so
\[
 |A_p(\varepsilon,\lambda)|
 \ll \min\{\varepsilon\log p,1\}\,p^{-1-\varepsilon},
\]
uniformly in $p$.  Hence $\sup_p|A_p|\to0$ and
$|A_p|^2\ll p^{-2}$.  Since $A_p\to0$ for each fixed prime, dominated
convergence gives $\sum_p|A_p|^2\to0$.

Let $K\subset\C$ be compact.  For small $\varepsilon$ the Taylor estimate
\[
 \log(1+zA)-z\log(1+A)=O_K(|A|^2)
\]
holds uniformly for $z\in K$ and $|A|\le\sup_p|A_p|$.  Thus
\[
 \sup_{z\in K}\left|
 \sum_p\bigl(\log(1+zA_p)-z\log(1+A_p)\bigr)\right|
 \ll_K\sum_p|A_p|^2\longrightarrow0.
\]
This proves that the product on the right of
\eqref{eq:marked-factorization-detail} tends to $1$ locally uniformly in
$z$.  On the positive real path,
\[
 \frac{\zeta(1+(1+2\lambda)\varepsilon)}{\zeta(1+\varepsilon)}
 \longrightarrow\frac1{1+2\lambda},
\]
so
\[
 \mathscr E_\ell(1+\varepsilon,s_\ell+\lambda\varepsilon;z)
 \longrightarrow(1+2\lambda)^{-z}
\]
locally uniformly in $z$.  For fixed $q$, choose a circle $|z|=R$ around
$0$; Cauchy's coefficient formula and uniform convergence on that circle give
\[
 [z^q]\mathscr E_\ell
 \longrightarrow[z^q](1+2\lambda)^{-z}
 =\frac{[-\log(1+2\lambda)]^q}{q!}.
\]
Finally,
$\varepsilon^r\widehat H(s_\ell+\lambda\varepsilon)
\to c_{-r}\lambda^{-r}$ by \eqref{eq:abstract-H-Laurent}, proving
\eqref{eq:abstract-fixedq}.  Expanding
$-\log(1+2\lambda)=-2\lambda+O(\lambda^2)$ gives
\eqref{eq:abstract-local-global}.
\end{proof}

The two results above describe complementary fixed-support aspects of the
same mechanism.  Theorem~\ref{thm:AMR-local} identifies locally which
prime-support layer survives, while Theorem~\ref{thm:AMR-fixedq} shows that
the same exponent persists after global summation when $q$ is fixed.  We now
change regime: $q$ is allowed to grow as the Dirichlet boundary is
approached, and the surviving jet weight itself becomes the object of study.

\section{Moving prime-support statistics}
\label{sec:prime-support-statistics}

In contrast with the fixed-$q$ regime of Section~\ref{sec:abstract-resonance},
the support depth now moves with $\varepsilon$.  We first derive the canonical
Euler-product law on the scale $q\varepsilon=O(1)$, including a uniform
$O(\varepsilon)$ coefficient estimate, and then normalize the same weights
probabilistically to obtain a mod-Poisson law with local asymptotics and a
Berry--Esseen bound.

\subsection{Moving prime support and a canonical Euler product}

The first surviving coefficient in Theorem~\ref{thm:AMR-local} singles out
the multiplicative weight
\begin{equation}\label{eq:jet-weight}
 W(N):=\prod_{p\mid N}\log p,
 \qquad W(1):=1.
\end{equation}
For $\Re\sigma>1$ put
\begin{equation}\label{eq:jet-weight-Euler}
 \mathcal F(\sigma,z)
 :=\sum_{N\ge1}\frac{W(N)z^{\omega(N)}}{N^\sigma}
 =\prod_p\left(1+z\frac{\log p}{p^\sigma-1}\right),
\end{equation}
and write
\[
 F_q(\sigma):=[z^q]\mathcal F(\sigma,z)
 =\sum_{\omega(N)=q}\frac{W(N)}{N^\sigma}.
\]
For $\varepsilon>0$ let
\[
 b_{p,\varepsilon}:=\frac{\log p}{p^{1+\varepsilon}-1}.
\]
For every fixed $\varepsilon>0$ one has
$\sum_p b_{p,\varepsilon}<\infty$, so
$\mathcal F(1+\varepsilon,\cdot)$ is entire.  Since
\begin{equation}\label{eq:logder-zeta-weight}
 \sum_p b_{p,\varepsilon}
 =-\frac{\zeta'(1+\varepsilon)}{\zeta(1+\varepsilon)}
 =\frac1\varepsilon-\gamma+O(\varepsilon),
\end{equation}
define
\begin{equation}\label{eq:Cepsilon}
 \mathcal C_\varepsilon(z)
 :=e^{-z/\varepsilon}\mathcal F(1+\varepsilon,z).
\end{equation}

The following theorem gives the global counterpart of the local resonance
law and identifies the canonical moving-support scale.

\begin{theorem}
\label{thm:moving-support}
The functions $\mathcal C_\varepsilon$ converge locally uniformly on $\C$ as
$\varepsilon\downarrow0$ to the entire function
\begin{equation}\label{eq:C0-product}
 \mathcal C_0(z)
 =e^{-\gamma z}\prod_p
 \left(1+\frac{z\log p}{p-1}\right)
 \exp\!\left(-\frac{z\log p}{p-1}\right).
\end{equation}
More precisely, for every compact set $K\subset\C$,
\begin{equation}\label{eq:Ceps-rate}
 \sup_{z\in K}|\mathcal C_\varepsilon(z)-\mathcal C_0(z)|
 =O_K(\varepsilon).
\end{equation}
Moreover, for every $R>0$,
\begin{equation}\label{eq:moving-support-uniform}
 \sup_{\substack{q\ge0\\q\varepsilon\le R}}
 \left|q!\varepsilon^qF_q(1+\varepsilon)
 -\mathcal C_0(q\varepsilon)\right|
 =O_R(\varepsilon).
\end{equation}
Consequently, if $q=q(\varepsilon)$ and
$q\varepsilon\to\rho\in[0,\infty)$, then
\begin{equation}\label{eq:moving-support-pointwise}
 q!\varepsilon^q
 \sum_{\omega(N)=q}\frac{W(N)}{N^{1+\varepsilon}}
 \longrightarrow\mathcal C_0(\rho).
\end{equation}
\end{theorem}

\begin{proof}
From \eqref{eq:logder-zeta-weight},
\begin{equation}\label{eq:Ceps-product}
 \mathcal C_\varepsilon(z)
 =\exp\!\left(
 z\left[-\frac{\zeta'(1+\varepsilon)}{\zeta(1+\varepsilon)}
 -\frac1\varepsilon\right]\right)
 \prod_p(1+zb_{p,\varepsilon})e^{-zb_{p,\varepsilon}}.
\end{equation}  Since
$0<b_{p,\varepsilon}\le \log p/(p-1)$ and
\[
 \sum_p\left(\frac{\log p}{p-1}\right)^2<\infty,
\]
the canonical products
$\prod_p(1+zb_{p,\varepsilon})e^{-zb_{p,\varepsilon}}$ converge locally
uniformly, uniformly for small $\varepsilon$.  We record the rate because it
is useful below.  Differentiation gives
\[
 \left|\frac{\partial}{\partial\varepsilon}b_{p,\varepsilon}\right|
 =\frac{(\log p)^2p^{1+\varepsilon}}
 {(p^{1+\varepsilon}-1)^2}
 \ll \frac{(\log p)^2}{p}
 \qquad(0\le\varepsilon\le\varepsilon_0).
\]
For $g_z(b):=(1+zb)e^{-zb}$ one has
$\partial_b g_z(b)=-z^2be^{-zb}$.  Hence, after separating finitely many
small primes, uniformly for $z$ in a compact set $K$,
\[
 \sum_p|g_z(b_{p,\varepsilon})-g_z(b_{p,0})|
 \ll_K \varepsilon
 \sum_p\frac{(\log p)^3}{p^2}
 =O_K(\varepsilon).
\]
The canonical products are locally uniformly bounded, so a telescoping
product estimate gives an $O_K(\varepsilon)$ difference between the two
products.  Together with
\[
 -\frac{\zeta'(1+\varepsilon)}{\zeta(1+\varepsilon)}
 -\frac1\varepsilon=-\gamma+O(\varepsilon),
\]
this proves \eqref{eq:C0-product} and the quantitative estimate
\eqref{eq:Ceps-rate}.

Write
$\mathcal C_\varepsilon(z)=\sum_{j\ge0}a_j(\varepsilon)z^j$.
Because
$\mathcal F(1+\varepsilon,z)=e^{z/\varepsilon}\mathcal C_\varepsilon(z)$,
coefficient extraction gives the exact identity
\begin{equation}\label{eq:coeff-identity-moving}
 q!\varepsilon^qF_q(1+\varepsilon)
 =\sum_{j=0}^q a_j(\varepsilon)(q)_j\varepsilon^j,
\end{equation}
where $(q)_j=q(q-1)\cdots(q-j+1)$.  Fix $R>0$ and choose $R_1>\max\{R,1\}$ sufficiently larger than $R$.
The estimate \eqref{eq:Ceps-rate} on $|z|\le R_1$, together with Cauchy's
formula, gives
\[
 |a_j(\varepsilon)-a_j(0)|\ll_{R_1}\varepsilon R_1^{-j},
 \qquad
 |a_j(\varepsilon)|+|a_j(0)|\ll_{R_1}R_1^{-j}.
\]
Put $x=q\varepsilon\le R$.  If $q\ge j$, then
\[
 0\le x^j-(q)_j\varepsilon^j
 \le \frac{j(j-1)}2\,\varepsilon x^{j-1};
\]
if $q<j$, the falling factorial vanishes and
$x^j\le j\varepsilon R^{j-1}$.  Thus in all cases
\[
 |(q)_j\varepsilon^j-x^j|
 \ll \varepsilon j^2\max\{1,R\}^{j-1}.
\]
Substituting these estimates in \eqref{eq:coeff-identity-moving} and summing
the resulting convergent geometric series proves the uniform
$O_R(\varepsilon)$ estimate \eqref{eq:moving-support-uniform}.  The pointwise
statement follows immediately.
\end{proof}

\begin{proposition}\label{prop:not-SD-admissible}
The multiplicative weight $W$ is not admissible, in the finite-average-value
Selberg--Delange sense used in \cite{Janisch2025} and developed in
\cite{dBT2021,Tenenbaum2015,Delange1954}, for any $\rho\in\C$.
More generally, for every fixed real $y>0$, the multiplicative weight
$N\mapsto y^{\omega(N)}W(N)$ is not admissible with any finite average
value.
\end{proposition}

\begin{proof}
The corresponding Dirichlet series is $\mathcal F(\sigma,y)$.  By
Theorem~\ref{thm:moving-support} and \eqref{eq:Ceps-rate},
\[
 \mathcal F(1+\varepsilon,y)
 =e^{y/\varepsilon}
  \bigl(\mathcal C_0(y)+O_y(\varepsilon)\bigr).
\]
For $y>0$ every factor in \eqref{eq:C0-product} is positive, so
$\mathcal C_0(y)>0$.  If this weight were Selberg--Delange admissible with
some finite average value $\rho$, then by definition
$\mathcal F(\sigma,y)\zeta(\sigma)^{-\rho}$ would extend holomorphically
across $\sigma=1$ and hence remain bounded on a sufficiently short real
interval to the right of $1$.  But along the positive real axis, interpreting
$\zeta(1+\varepsilon)^{-\rho}$ via the real logarithm of the positive number
$\zeta(1+\varepsilon)$, one has
$|\zeta(1+\varepsilon)^{-\rho}|\sim\varepsilon^{\operatorname{Re}\rho}$,
while the displayed formula grows like $e^{y/\varepsilon}$.  Thus the
product is unbounded as
$\varepsilon\downarrow0$, a contradiction.  Taking $y=1$ gives the claim
for $W$.
\end{proof}

The moving law also applies directly to the first surviving jets of every
kernel covered by the local transfer theorem.

\begin{corollary}
\label{cor:moving-first-jets}
Under the hypotheses of Theorem~\ref{thm:AMR-local}, fix $r$ and put, for
$q\ge r$,
\[
 D_{q,r,H}(N)
 :=\left.\frac{\dd^{q-r}}{\dd s^{q-r}}
 \widehat{\mathscr P_{N,\ell}H}(s)\right|_{s=s_\ell}
 \qquad(\omega(N)=q).
\]
Then, for every $R>0$,
\begin{equation}\label{eq:moving-first-jets}
 \sup_{\substack{q\ge r\\q\varepsilon\le R}}
 \left|
 \left(\frac{q!}{(q-r)!}\right)
 \left(\frac{(-\varepsilon/2)^q}{c_{-r}}\right)
 \sum_{\omega(N)=q}\frac{D_{q,r,H}(N)}{N^{1+\varepsilon}}
 -\mathcal C_0(q\varepsilon)
 \right|=O_R(\varepsilon).
\end{equation}
Thus, after the explicit jet normalization in
\eqref{eq:moving-first-jets}, the moving prime-support law is independent of
the kernel $H$.
\end{corollary}

\begin{proof}
By \eqref{eq:abstract-first-jet} and \eqref{eq:abstract-support},
\[
 D_{q,r,H}(N)=(q-r)!(-2)^q c_{-r}W(N).
\]
The left side of \eqref{eq:moving-first-jets} is therefore exactly the
quantity in \eqref{eq:moving-support-uniform}.
\end{proof}

\subsection{A mod-Poisson law forced by Mellin jets}

The probability law in this subsection is intrinsic to the resonance
mechanism: its weights are precisely the first nonzero Jordan--Mellin
coefficients on the layers $\omega(N)=q$.  Normalizing these coefficients
therefore turns the arithmetic filtration itself into a probabilistic
ensemble.  The next theorem identifies the limiting law together with an
effective normal approximation and local asymptotics.

For $\varepsilon>0$ define a random positive integer
$\mathsf N_\varepsilon$ by
\begin{equation}\label{eq:Dirichlet-jet-probability}
 \Prob(\mathsf N_\varepsilon=N)
 :=\frac{W(N)N^{-1-\varepsilon}}
 {\mathcal F(1+\varepsilon,1)},
 \qquad N\ge1,
\end{equation}
and put $\Omega_\varepsilon:=\omega(\mathsf N_\varepsilon)$.

\begin{theorem}
\label{thm:mod-Poisson-jet}
As $\varepsilon\downarrow0$,
\begin{equation}\label{eq:mod-Poisson-pgf}
 e^{-(z-1)/\varepsilon}\E[z^{\Omega_\varepsilon}]
 =\Psi(z)+O_K(\varepsilon),
 \qquad
 \Psi(z):=\frac{\mathcal C_0(z)}{\mathcal C_0(1)},
\end{equation}
uniformly for $z$ in every compact set $K\subset\C$.  Equivalently,
uniformly for $u$ in compact subsets of $\C$,
\begin{equation}\label{eq:mod-Poisson-mgf}
 \E[e^{u\Omega_\varepsilon}]
 =\exp\!\left(\frac{e^u-1}{\varepsilon}\right)
 \bigl(\Psi(e^u)+O_K(\varepsilon)\bigr).
\end{equation}
Hence $\Omega_\varepsilon$ has mod-Poisson convergence with parameter
$1/\varepsilon$.  In particular,
\begin{equation}\label{eq:jet-CLT}
 \sqrt\varepsilon\left(
 \Omega_\varepsilon-\frac1\varepsilon\right)
 \xrightarrow{\ d\ }\mathcal N(0,1).
\end{equation}
In fact, if $d_{\rm K}$ denotes Kolmogorov distance, then
\begin{equation}\label{eq:jet-Berry-Esseen}
 d_{\rm K}\!\left(
 \frac{\Omega_\varepsilon-1/\varepsilon}{\sqrt{1/\varepsilon}},
 \mathcal N(0,1)\right)=O(\sqrt\varepsilon).
\end{equation}
More precisely, for every $R>0$, uniformly for $q\varepsilon\le R$,
\begin{equation}\label{eq:jet-local-deviation}
 \Prob(\Omega_\varepsilon=q)
 =
 \frac{\mathcal C_0(q\varepsilon)}{\mathcal C_0(1)}
 e^{-1/\varepsilon}\frac{\varepsilon^{-q}}{q!}
 \bigl(1+O_R(\varepsilon)\bigr).
\end{equation}
In particular, if $q\varepsilon\to\rho\in[0,\infty)$, this gives the
corresponding precise local asymptotic with $\mathcal C_0(\rho)$.
For $\rho>0$, Stirling's formula yields the local large-deviation rate
\begin{equation}\label{eq:jet-LDP-rate}
 \varepsilon\log\Prob(\Omega_\varepsilon=q)
 \longrightarrow-[\rho\log\rho-\rho+1].
\end{equation}
\end{theorem}

\begin{proof}
By \eqref{eq:Cepsilon},
\[
 \E[z^{\Omega_\varepsilon}]
 =\frac{\mathcal F(1+\varepsilon,z)}
 {\mathcal F(1+\varepsilon,1)}
 =e^{(z-1)/\varepsilon}
 \frac{\mathcal C_\varepsilon(z)}{\mathcal C_\varepsilon(1)},
\]
and \eqref{eq:Ceps-rate} proves the quantitative forms
\eqref{eq:mod-Poisson-pgf} and \eqref{eq:mod-Poisson-mgf}.  For the central limit theorem, put
$z=e^{t\sqrt\varepsilon}$ and use
\[
 \frac{e^{t\sqrt\varepsilon}-1-t\sqrt\varepsilon}{\varepsilon}
 \longrightarrow\frac{t^2}{2},
\]
while
$\mathcal C_\varepsilon(e^{t\sqrt\varepsilon})/
 \mathcal C_\varepsilon(1)\to1$.  The moment-generating functions therefore
converge to $e^{t^2/2}$ near the origin.

For the effective normal approximation, factor the probability generating
function prime by prime:
\[
 \E[z^{\Omega_\varepsilon}]
 =\prod_p\frac{1+zb_{p,\varepsilon}}{1+b_{p,\varepsilon}}
 =\prod_p\bigl((1-\pi_{p,\varepsilon})+
                 \pi_{p,\varepsilon}z\bigr),
 \qquad
 \pi_{p,\varepsilon}:=\frac{b_{p,\varepsilon}}{1+b_{p,\varepsilon}}.
\]
Thus $\Omega_\varepsilon$ has the law of a convergent sum of independent
Bernoulli variables.  Since $\sum_pb_{p,\varepsilon}^2=O(1)$,
\[
 \mu_\varepsilon:=\sum_p\pi_{p,\varepsilon}
 =\frac1\varepsilon+O(1),\qquad
 \sigma_\varepsilon^2:=\sum_p\pi_{p,\varepsilon}
 (1-\pi_{p,\varepsilon})=\frac1\varepsilon+O(1).
\]
The sum of the third absolute centered moments is $O(1/\varepsilon)$.
Applying the classical Berry--Esseen theorem to finite prime truncations
and then passing to the $L^2$ limit gives
\[
 d_{\rm K}\!\left(
 \frac{\Omega_\varepsilon-\mu_\varepsilon}{\sigma_\varepsilon},
 \mathcal N(0,1)\right)=O(\sqrt\varepsilon).
\]
Since $\mu_\varepsilon=1/\varepsilon+O(1)$ and
$\sigma_\varepsilon=\varepsilon^{-1/2}(1+O(\varepsilon))$, the deterministic
change of centering and scale contributes only $O(\sqrt\varepsilon)$ to the
Kolmogorov distance.  This proves \eqref{eq:jet-Berry-Esseen}.

Finally,
\[
 \Prob(\Omega_\varepsilon=q)
 =\frac{F_q(1+\varepsilon)}{\mathcal F(1+\varepsilon,1)}.
\]
Theorem~\ref{thm:moving-support} gives, uniformly for
$q\varepsilon\le R$,
\[
 F_q(1+\varepsilon)
 =\frac{\varepsilon^{-q}}{q!}
  \bigl(\mathcal C_0(q\varepsilon)+O_R(\varepsilon)\bigr),
\]
whereas \eqref{eq:Ceps-rate} at $z=1$ gives
\[
 \mathcal F(1+\varepsilon,1)
 =e^{1/\varepsilon}\mathcal C_0(1)
  \bigl(1+O(\varepsilon)\bigr).
\]
Since $\mathcal C_0$ is positive and bounded away from zero on every
compact interval $[0,R]$, division proves \eqref{eq:jet-local-deviation};
Stirling's formula proves \eqref{eq:jet-LDP-rate}.
\end{proof}

\section{Character twists, Frobenius selection, and Galois entanglement}
\label{sec:twisted-resonance}

Sections~\ref{sec:abstract-resonance} and~\ref{sec:prime-support-statistics} count all distinct prime divisors.  The
resonance mechanism admits an arithmetic refinement that distinguishes
classes of primes.  After twisting by a Dirichlet character, the critical
vanishing order is determined only by the selected primes satisfying
$\chi(p)=1$.  Let $\chi$ be a primitive Dirichlet character modulo $Q$.  To avoid the genuinely
degenerate local factors at ramified primes, throughout the local statements
below we assume $(N,Q)=1$.

For $\ell\in\mathbb N$ define
\begin{equation}\label{eq:twisted-P}
 (\mathscr P_{N,\ell,\chi}H)(t)
 :=\sum_{d\mid N}\mu_{\rm M}(N/d)\chi(d)d^\ell H(d^2t),
\end{equation}
and put
\begin{equation}\label{eq:twisted-Jordan}
 J_{\nu,\chi}(N)
 :=\sum_{d\mid N}\mu_{\rm M}(N/d)\chi(d)d^\nu,
 \qquad
 \Lambda_{k,\chi}(N)
 :=\sum_{d\mid N}\mu_{\rm M}(N/d)\chi(d)(\log d)^k.
\end{equation}
Let
\begin{align*}
 S_\chi(N)&:=\{p\mid N:\chi(p)=1\},\qquad
 \omega_{\chi,1}(N):=\#S_\chi(N),\\
 W_\chi(N)&:=\prod_{p\in S_\chi(N)}\log p.
\end{align*}
and define the nonvanishing complementary factor
\begin{equation}\label{eq:Cchi-local}
 C_\chi(N)
 :=\prod_{p^a\Vert N,\,\chi(p)\ne1}
   \chi(p)^{a-1}(\chi(p)-1),
\end{equation}
with the empty product equal to one.

\begin{theorem}\label{thm:twisted-local-resonance}
Under the hypotheses of Theorem~\ref{thm:AMR-local}, and for $(N,Q)=1$,
\begin{equation}\label{eq:twisted-Mellin-factor}
 \widehat{\mathscr P_{N,\ell,\chi}H}(s)
 =J_{\ell-2s,\chi}(N)\widehat H(s).
\end{equation}
At $s=s_\ell+u$,
\begin{equation}\label{eq:twisted-Jordan-jet}
 J_{-2u,\chi}(N)
 =\sum_{k\ge0}\frac{(-2)^k}{k!}
   \Lambda_{k,\chi}(N)u^k.
\end{equation}
If $q_\chi=\omega_{\chi,1}(N)$, then
\begin{equation}\label{eq:twisted-first-Lambda}
 \Lambda_{k,\chi}(N)=0\quad(k<q_\chi),\qquad
 \Lambda_{q_\chi,\chi}(N)
 =q_\chi!\,C_\chi(N)W_\chi(N),
\end{equation}
and consequently
\begin{equation}\label{eq:twisted-valuation}
 \ord_{s=s_\ell}
 \widehat{\mathscr P_{N,\ell,\chi}H}
 =\omega_{\chi,1}(N)-r.
\end{equation}
If $q_\chi\ge r$, the first surviving derivative is
\begin{equation}\label{eq:twisted-first-jet}
 \left.\frac{\dd^{q_\chi-r}}{\dd s^{q_\chi-r}}
 \widehat{\mathscr P_{N,\ell,\chi}H}(s)\right|_{s=s_\ell}
 =(q_\chi-r)!(-2)^{q_\chi}
 C_\chi(N)W_\chi(N)c_{-r}.
\end{equation}
\end{theorem}

\begin{proof}
The Mellin scaling argument is the same as in
Theorem~\ref{thm:AMR-local}.  Multiplicativity gives, for $p^a\Vert N$,
\begin{equation}\label{eq:twisted-prime-power}
 J_{\nu,\chi}(p^a)
 =\chi(p)^{a-1}p^{(a-1)\nu}
   \bigl(\chi(p)p^\nu-1\bigr).
\end{equation}
When $\nu=-2u$, a prime with $\chi(p)=1$ contributes
$-2u\log p+O(u^2)$, whereas a prime with $\chi(p)\ne1$ contributes the
nonzero constant $\chi(p)^{a-1}(\chi(p)-1)$ at $u=0$.  Hence
\[
 J_{-2u,\chi}(N)
 =(-2u)^{q_\chi}C_\chi(N)W_\chi(N)
  +O_N(u^{q_\chi+1}).
\]
Comparison with the Taylor expansion
\eqref{eq:twisted-Jordan-jet} proves
\eqref{eq:twisted-first-Lambda}.  Multiplication by the Laurent expansion of
$\widehat H$ then gives \eqref{eq:twisted-valuation} and
\eqref{eq:twisted-first-jet}.
\end{proof}

\begin{remark}\label{rem:twisted-ramified}
The condition $(N,Q)=1$ is structural rather than cosmetic.  If $p\mid Q$,
then $\chi(p)=0$; the local factor equals $-1$ when $p\Vert N$ but vanishes
identically when $p^2\mid N$.  Thus ramified prime powers create a different
degeneracy from the resonance zeros studied here.
\end{remark}

\begin{corollary}\label{cor:quadratic-splitting}
Let $D$ be a fundamental discriminant and let $\chi_D$ be the associated
primitive quadratic character.  For $(N,D)=1$,
\[
 \omega_{\chi_D,1}(N)
 =\#\{p\mid N:p\text{ splits in }\mathbb Q(\sqrt D)\}.
\]
Hence the twisted resonance filtration counts split prime divisors and ignores
inert prime divisors at the level of the vanishing order.  More generally, for
a Dirichlet character $\chi$, the condition $\chi(p)=1$ means that the
Frobenius class of $p$ in the cyclotomic extension lies in $\ker\chi$.
\end{corollary}

\begin{proof}
For an unramified prime $p$, the Kronecker symbol satisfies
$\chi_D(p)=1$ exactly when $p$ splits in the quadratic field, and equals
$-1$ when $p$ is inert.  The final statement is the standard interpretation
of a Dirichlet character as a one-dimensional character of the abelian
Galois group of a cyclotomic extension; see \cite{Davenport2000}.
\end{proof}

The local Frobenius selection has a global probabilistic counterpart.
Restricting to the selected-prime subsemigroup yields a positive ensemble
whose natural scale is determined by the density $1/h$ of the condition
$\chi(p)=1$.  Assume
that $\chi$ is nonprincipal of exact order $h\ge2$, and let
\[
 \mathcal P_\chi:=\{p:p\nmid Q,\ \chi(p)=1\},\qquad
 \delta_\chi:=\frac1h.
\]
For integers all of whose prime divisors belong to $\mathcal P_\chi$, put
$W_\chi(N)=\prod_{p\mid N}\log p$, and define
\begin{equation}\label{eq:Fchi-Euler}
 \mathcal F_\chi(\sigma,z)
 :=\prod_{p\in\mathcal P_\chi}
 \left(1+z\frac{\log p}{p^\sigma-1}\right),
 \qquad
 F_{\chi,q}(\sigma):=[z^q]\mathcal F_\chi(\sigma,z).
\end{equation}

\begin{theorem}\label{thm:Frobenius-moving-support}
There exists a real constant $\kappa_\chi$ such that
\begin{equation}\label{eq:split-prime-sum}
 \sum_{p\in\mathcal P_\chi}
 \frac{\log p}{p^{1+\varepsilon}-1}
 =\frac{\delta_\chi}{\varepsilon}+\kappa_\chi+O_\chi(\varepsilon).
\end{equation}
Consequently
\begin{equation}\label{eq:Cchi-eps}
 \mathcal C_{\chi,\varepsilon}(z)
 :=e^{-\delta_\chi z/\varepsilon}
   \mathcal F_\chi(1+\varepsilon,z)
\end{equation}
converges locally uniformly, with $O_{\chi,K}(\varepsilon)$ error on compact
sets, to the entire function
\begin{equation}\label{eq:Cchi-zero}
 \mathcal C_{\chi,0}(z)
 =e^{\kappa_\chi z}
  \prod_{p\in\mathcal P_\chi}
  \left(1+\frac{z\log p}{p-1}\right)
  \exp\!\left(-\frac{z\log p}{p-1}\right).
\end{equation}
For every $R>0$,
\begin{equation}\label{eq:Frobenius-moving-law}
 \sup_{q\varepsilon/\delta_\chi\le R}
 \left|
 q!\left(\frac{\varepsilon}{\delta_\chi}\right)^q
 F_{\chi,q}(1+\varepsilon)
 -\mathcal C_{\chi,0}\!\left(\frac{q\varepsilon}{\delta_\chi}\right)
 \right|
 =O_{\chi,R}(\varepsilon).
\end{equation}
If $\mathsf N_{\chi,\varepsilon}$ is distributed on the split-supported
integers by
\[
 \Prob(\mathsf N_{\chi,\varepsilon}=N)
 =\frac{W_\chi(N)N^{-1-\varepsilon}}
 {\mathcal F_\chi(1+\varepsilon,1)},
\]
and $\Omega_{\chi,\varepsilon}=\omega(\mathsf N_{\chi,\varepsilon})$, then
\begin{equation}\label{eq:Frobenius-mod-Poisson}
 e^{-\delta_\chi(z-1)/\varepsilon}
 \E[z^{\Omega_{\chi,\varepsilon}}]
 \longrightarrow
 \frac{\mathcal C_{\chi,0}(z)}{\mathcal C_{\chi,0}(1)}
\end{equation}
locally uniformly on $\C$.  Thus the mod-Poisson parameter is
$\delta_\chi/\varepsilon=1/(h\varepsilon)$; for a quadratic character it is
$1/(2\varepsilon)$.
\end{theorem}

\begin{proof}
For $\Re s>1$, first separate the $a=1$ term from the prime-power tail:
\[
 \sum_{p\in\mathcal P_\chi}\frac{\log p}{p^s-1}
 =\sum_{p\in\mathcal P_\chi}\frac{\log p}{p^s}
  +H_\chi(s),
\]
where $H_\chi$ is holomorphic in a neighborhood of $s=1$.  Character
orthogonality gives, for $p\nmid Q$,
\[
 \mathbf1_{\{\chi(p)=1\}}
 =\frac1h\sum_{j=0}^{h-1}\chi(p)^j.
\]
Let $\chi^0$ denote the principal character modulo $Q$.  For
$0\le j<h$, the Euler product gives
\[
 -\frac{L'(s,\chi^j)}{L(s,\chi^j)}
 =\sum_{p\nmid Q}\sum_{a\ge1}
   \frac{\chi(p)^{ja}\log p}{p^{as}}.
\]
Hence, for $\Re s>1$,
\begin{equation}\label{eq:Frobenius-Llog-explicit}
 \sum_{p\in\mathcal P_\chi}\frac{\log p}{p^s}
 =-\frac1h\sum_{j=0}^{h-1}
   \frac{L'(s,\chi^j)}{L(s,\chi^j)}+H_\chi^*(s),
\end{equation}
where
\[
 H_\chi^*(s)
 :=-\frac1h\sum_{j=0}^{h-1}\sum_{p\nmid Q}\sum_{a\ge2}
   \frac{\chi(p)^{ja}\log p}{p^{as}}
\]
is holomorphic for $\Re s>1/2$.  Since
$L(s,\chi^0)=\zeta(s)\prod_{p\mid Q}(1-p^{-s})$, the term $j=0$ contributes
$1/(s-1)$ at $s=1$.  For $1\le j<h$, exact order $h$ implies that $\chi^j$
is nonprincipal, so $L(s,\chi^j)$ is holomorphic and nonzero at $s=1$.
Combining \eqref{eq:Frobenius-Llog-explicit} with the holomorphic prime-power
tail $H_\chi(s)$ gives a Laurent expansion
$\delta_\chi/(s-1)+\kappa_\chi+O_\chi(s-1)$, with
$\delta_\chi=1/h$.  Setting $s=1+\varepsilon$ proves
\eqref{eq:split-prime-sum}; see \cite{Davenport2000} for the standard
$L$-function facts.

The canonical-product proof of Theorem~\ref{thm:moving-support} now applies
prime-for-prime to the subset $\mathcal P_\chi$, giving
\eqref{eq:Cchi-zero} and the local uniform $O(\varepsilon)$ rate.  Coefficient
extraction from
$\mathcal F_\chi=e^{\delta_\chi z/\varepsilon}
\mathcal C_{\chi,\varepsilon}$ gives
\eqref{eq:Frobenius-moving-law} after the normalization
$\varepsilon/\delta_\chi$.  Finally, division by the value at $z=1$ gives
\eqref{eq:Frobenius-mod-Poisson}.
\end{proof}

\begin{remark}
The local theorem allows inert or otherwise nonselected prime divisors to be
present: they contribute the nonzero factor $C_\chi(N)$ but do not increase
the resonance depth.  The positive global ensemble in
Theorem~\ref{thm:Frobenius-moving-support} restricts to the selected-prime
subsemigroup so that the first surviving weights are positive.  Thus the two
statements separate the local Galois selection mechanism from its positive
probabilistic realization.
\end{remark}

\subsection{Joint Frobenius support and Galois entanglement}

The preceding one-dimensional selection becomes coupled in the splitting-count
coordinates when the same prime divisor is allowed to satisfy several splitting
conditions at once.
Let $L_1,\ldots,L_r/\mathbb Q$ be finite Galois extensions, let
$L=L_1\cdots L_r$, and write $G=\operatorname{Gal}(L/\mathbb Q)$.  We omit the
finitely many primes dividing the discriminant $\Delta_L$.  For
$p\nmid\Delta_L$ put
\[
 m_i(p):=\mathbf1_{\{p\text{ splits completely in }L_i\}},
 \qquad
 \mathbf m(p):=(m_1(p),\ldots,m_r(p)).
\]
For $g\in G$ define $m_i(g)=1$ if the restriction of $g$ to $L_i$ is the
identity and $m_i(g)=0$ otherwise.  This is a class function and
$\mathbf m(p)=\mathbf m(\operatorname{Frob}_p)$.

We first record the fixed-extension analytic input.

\begin{lemma}
\label{lem:Dedekind-splitting-prime-sum}
Let $K/\mathbb Q$ be a finite Galois extension of degree $n_K$ and
discriminant $\Delta_K$.  There is a constant $\kappa_K$ such that, as
$\varepsilon\downarrow0$,
\begin{equation}\label{eq:Dedekind-splitting-prime-sum}
 \sum_{\substack{p\nmid\Delta_K\\
                  p\text{ splits completely in }K}}
 \frac{\log p}{p^{1+\varepsilon}-1}
 =\frac1{n_K\varepsilon}+\kappa_K+O_K(\varepsilon).
\end{equation}
The same statement holds for $K=\mathbb Q$, with the splitting condition
vacuous and $n_K=1$.
\end{lemma}

\begin{proof}
For $\Re s>1$, the logarithmic derivative of the Dedekind zeta function is
\[
 -\frac{\zeta_K'}{\zeta_K}(s)
 =\sum_{\mathfrak p}\sum_{m\ge1}
   \frac{\log N\mathfrak p}{(N\mathfrak p)^{ms}}.
\]
If $p\nmid\Delta_K$ is unramified, the Galois hypothesis implies that all
prime ideals above $p$ have the same residue degree $f_p$.  Their total
contribution is therefore
\[
 n_K\frac{\log p}{p^{f_ps}-1}.
\]
For primes that split completely, $f_p=1$.  For every other unramified prime,
$f_p\ge2$.  Hence
\begin{equation}\label{eq:Dedekind-log-derivative-split}
 -\frac1{n_K}\frac{\zeta_K'}{\zeta_K}(s)
 =\sum_{\substack{p\nmid\Delta_K\\
                  p\text{ splits completely in }K}}
   \frac{\log p}{p^s-1}+H_K(s),
\end{equation}
where $H_K$ consists of the non-split unramified primes together with the
finitely many ramified local factors.  The non-split part converges normally
for $\Re s>1/2$, since
\[
 \sum_p\frac{\log p}{p^{2\Re s}-1}
\]
converges there; thus $H_K$ is holomorphic near $s=1$.

The Dedekind zeta function has a simple pole and no zero at $s=1$, so
\[
 -\frac{\zeta_K'}{\zeta_K}(s)
 =\frac1{s-1}+c_K+O_K(s-1).
\]
Substitution into \eqref{eq:Dedekind-log-derivative-split} and
$s=1+\varepsilon$ proves \eqref{eq:Dedekind-splitting-prime-sum}.  These are
standard properties of Dedekind zeta functions; see, for example,
\cite[Chap.~VII]{Neukirch1999}.
\end{proof}

For $\mathbf z=(z_1,\ldots,z_r)$ write
$\mathbf z^{\mathbf m}=\prod_i z_i^{m_i}$ and define
\begin{equation}\label{eq:Galois-joint-Euler}
 \mathcal F_{\mathbf L}(\sigma,\mathbf z)
 :=\prod_{p\nmid\Delta_L}
 \left(1+\frac{\log p}{p^\sigma-1}
       \mathbf z^{\mathbf m(p)}\right).
\end{equation}
For a nonempty set $S\subseteq[r]:=\{1,\ldots,r\}$ let
\[
 L_S:=\prod_{i\in S}L_i,
 \qquad
 \delta_S:=\frac1{[L_S:\mathbb Q]},
\]
and put
\begin{equation}\label{eq:Galois-leading-polynomial}
 P_{\mathbf L}(\mathbf z)
 :=\frac1{|G|}\sum_{g\in G}\mathbf z^{\mathbf m(g)}
 =1+\sum_{\emptyset\ne S\subseteq[r]}
   \delta_S\prod_{i\in S}(z_i-1).
\end{equation}
The second identity follows by expanding
$\prod_i(1+m_i(g)(z_i-1))$: for each $S$, the average of
$\prod_{i\in S}m_i(g)$ is
$|\operatorname{Gal}(L/L_S)|/|G|=[L_S:\mathbb Q]^{-1}$.

\begin{remark}
For $\mathbf v\in\{0,1\}^r$ let
$\mathcal P_{\mathbf v}:=\{p\nmid\Delta_L:\mathbf m(p)=\mathbf v\}$.
These sets partition the unramified primes, and therefore
\[
 \mathcal F_{\mathbf L}(\sigma,\mathbf z)
 =\prod_{\mathbf v}
   \prod_{p\in\mathcal P_{\mathbf v}}
   \left(1+\frac{\log p}{p^\sigma-1}\mathbf z^{\mathbf v}\right).
\]
Accordingly, under the weighted boundary measure, the exact-pattern support
counts are independent; the vector $(M_1,\ldots,M_r)$ is their linear image
$\sum_{\mathbf v}\mathbf v K_{\mathbf v}$.  Thus the coupling studied below
is not a claim of primitive dependence between the disjoint Frobenius pattern
classes.  Its arithmetic content is that the pattern intensities and all mixed
coefficients of $P_{\mathbf L}$ are constrained by the common Galois
compositum, which is what produces the shared jumps, the support lattice, and
the compositum-degree cumulants.
\end{remark}

\begin{theorem}
\label{thm:Galois-entanglement-support}
There is an entire function $\mathcal C_{\mathbf L,0}$ on $\mathbb C^r$ such
that, for every compact $K\subset\mathbb C^r$,
\begin{equation}\label{eq:Galois-canonical-limit}
 \sup_{\mathbf z\in K}
 \left|
 e^{-P_{\mathbf L}(\mathbf z)/\varepsilon}
 \mathcal F_{\mathbf L}(1+\varepsilon,\mathbf z)
 -\mathcal C_{\mathbf L,0}(\mathbf z)
 \right|
 =O_{\mathbf L,K}(\varepsilon).
\end{equation}

Let $\mathsf N_\varepsilon$ be distributed on the positive integers coprime to
$\Delta_L$ by
\begin{equation}\label{eq:Galois-weighted-ensemble}
 \Prob(\mathsf N_\varepsilon=N)
 =\frac{W(N)N^{-1-\varepsilon}}
 {\mathcal F_{\mathbf L}(1+\varepsilon,\mathbf1)},
 \qquad
 W(N)=\prod_{p\mid N}\log p,
\end{equation}
and set
\[
 M_{i,\varepsilon}
 :=\#\{p\mid\mathsf N_\varepsilon:
          p\text{ splits completely in }L_i\}.
\]
Then, locally uniformly for $\mathbf z$ near $\mathbf1$,
\begin{equation}\label{eq:Galois-mod-compound-Poisson}
 e^{-[P_{\mathbf L}(\mathbf z)-1]/\varepsilon}
 \E\!\left[\prod_{i=1}^r z_i^{M_{i,\varepsilon}}\right]
 \longrightarrow
 \Psi_{\mathbf L}(\mathbf z)
 :=\frac{\mathcal C_{\mathbf L,0}(\mathbf z)}
 {\mathcal C_{\mathbf L,0}(\mathbf1)}.
\end{equation}
Equivalently, the renormalizing reference law is compound Poisson with L\'evy masses
$\pi_{\mathbf v}/\varepsilon$ on the nonzero jump vectors
$\mathbf v\in\{0,1\}^r\setminus\{0\}$, where
$\pi_{\mathbf v}=|G|^{-1}\#\{g:\mathbf m(g)=\mathbf v\}$.
Thus \eqref{eq:Galois-mod-compound-Poisson} is a mod-compound-Poisson
convergence: the reference intensity is of order $\varepsilon^{-1}$ and hence
does not converge to a fixed compound-Poisson distribution.

For a multi-index $\boldsymbol\alpha\in\mathbb Z_{\ge0}^r\setminus\{0\}$,
let $S(\boldsymbol\alpha)=\{i:\alpha_i>0\}$ and denote by
$\kappa_{\boldsymbol\alpha}(\mathbf M_\varepsilon)$ the corresponding joint
cumulant.  Then
\begin{equation}\label{eq:Galois-cumulant-compositum}
 \kappa_{\boldsymbol\alpha}(\mathbf M_\varepsilon)
 =\frac{1}{\varepsilon[L_{S(\boldsymbol\alpha)}:\mathbb Q]}+O_{\mathbf L,\boldsymbol\alpha}(1).
\end{equation}
In particular,
\begin{equation}\label{eq:Galois-means-covariance}
 \E M_{i,\varepsilon}
 =\frac{1}{\varepsilon[L_i:\mathbb Q]}+O(1),
 \qquad
 \operatorname{Cov}(M_{i,\varepsilon},M_{j,\varepsilon})
 =\frac{1}{\varepsilon[L_iL_j:\mathbb Q]}+O(1).
\end{equation}
Consequently, for $i\ne j$,
\begin{equation}\label{eq:Galois-intersection-recovery}
 \lim_{\varepsilon\downarrow0}
 \frac{\varepsilon\operatorname{Cov}(M_{i,\varepsilon},M_{j,\varepsilon})}
 {(\varepsilon\E M_{i,\varepsilon})(\varepsilon\E M_{j,\varepsilon})}
 =[L_i\cap L_j:\mathbb Q].
\end{equation}
More generally, for every nonempty $S\subseteq[r]$,
\begin{equation}\label{eq:Galois-entanglement-factor}
 \lim_{\varepsilon\downarrow0}
 \frac{\varepsilon\,\kappa(M_{i,\varepsilon}:i\in S)}
 {\prod_{i\in S}(\varepsilon\E M_{i,\varepsilon})}
 =\frac{\prod_{i\in S}[L_i:\mathbb Q]}{[L_S:\mathbb Q]}.
\end{equation}
\end{theorem}

\begin{proof}
For every unramified prime,
\[
 \mathbf z^{\mathbf m(p)}
 =\prod_{i=1}^r\bigl(1+m_i(p)(z_i-1)\bigr)
 =\sum_{S\subseteq[r]}
   \left(\prod_{i\in S}(z_i-1)\right)
   \mathbf1_{\{p\text{ splits completely in }L_S\}},
\]
where $L_\emptyset:=\mathbb Q$.  Indeed, simultaneous complete splitting in
the fields $L_i$, $i\in S$, is equivalent to complete splitting in their
compositum $L_S$.  Applying
Lemma~\ref{lem:Dedekind-splitting-prime-sum} to the finitely many fields
$L_S$ gives the required expansions.  Replacing the condition
$p\nmid\Delta_{L_S}$ by the common restriction $p\nmid\Delta_L$ removes only
finitely many primes and therefore changes only the constant term.  Hence,
uniformly for $\mathbf z$ in a compact set,
\begin{equation}\label{eq:Galois-first-prime-sum}
 \sum_{p\nmid\Delta_L}
 \frac{\mathbf z^{\mathbf m(p)}\log p}{p^{1+\varepsilon}-1}
 =\frac{P_{\mathbf L}(\mathbf z)}{\varepsilon}
  +\kappa_{\mathbf L}(\mathbf z)+O_{\mathbf L,K}(\varepsilon),
\end{equation}
where $\kappa_{\mathbf L}$ is a polynomial in $\mathbf z$.  This proves the
required prime-sum expansion without any appeal to Artin $L$-function
holomorphy.

Put $b_{p,\varepsilon}=\log p/(p^{1+\varepsilon}-1)$.  As in the proof of
Theorem~\ref{thm:moving-support}, factor
\begin{align*}
&e^{-P_{\mathbf L}(\mathbf z)/\varepsilon}
 \mathcal F_{\mathbf L}(1+\varepsilon,\mathbf z)\\
&\qquad=
 \exp\!\left(
  \sum_{p\nmid\Delta_L}b_{p,\varepsilon}\mathbf z^{\mathbf m(p)}
  -\frac{P_{\mathbf L}(\mathbf z)}{\varepsilon}
 \right)
 \prod_{p\nmid\Delta_L}
 \bigl(1+b_{p,\varepsilon}\mathbf z^{\mathbf m(p)}\bigr)
 e^{-b_{p,\varepsilon}\mathbf z^{\mathbf m(p)}}.
\end{align*}
The first factor converges with $O_K(\varepsilon)$ speed by
\eqref{eq:Galois-first-prime-sum}.  For the canonical product, normal
convergence follows from
$\sum_p(\log p/(p-1))^2<\infty$.  Moreover,
$|\partial_\varepsilon b_{p,\varepsilon}|\ll(\log p)^2/p$, and differentiation
of $(1+bw)e^{-bw}$ contributes an additional factor $b$; hence the difference
between the $\varepsilon$-product and the limiting product is bounded on
compact sets by
\[
 O_{\mathbf L,K}(\varepsilon)
 \sum_p\frac{(\log p)^3}{p^2}=O_{\mathbf L,K}(\varepsilon).
\]
This proves \eqref{eq:Galois-canonical-limit}.

The Euler product \eqref{eq:Galois-joint-Euler} is exactly the joint
probability-generating series of the variables $M_{i,\varepsilon}$ under
\eqref{eq:Galois-weighted-ensemble}.  At $\mathbf z=\mathbf1$ the limiting
canonical product is strictly positive, so it remains nonzero in a neighborhood
of $\mathbf1$.  Division by the value at $\mathbf z=\mathbf1$ therefore gives
\eqref{eq:Galois-mod-compound-Poisson}.  Since
$P_{\mathbf L}(\mathbf1)=1$ and
\[
 P_{\mathbf L}(\mathbf z)-1
 =\sum_{\mathbf v\in\{0,1\}^r\setminus\{0\}}
 \pi_{\mathbf v}(\mathbf z^{\mathbf v}-1),
 \qquad
 \pi_{\mathbf v}:=\frac1{|G|}\#\{g:\mathbf m(g)=\mathbf v\},
\]
the reference law is compound Poisson with the stated L\'evy masses.

Finally set $z_i=e^{t_i}$ and take logarithms in
\eqref{eq:Galois-mod-compound-Poisson}.  Shrinking the neighborhood if
necessary, the limiting factor is nonzero there; local uniform convergence and
Cauchy's formula therefore control all derivatives of its logarithm.  Hence,
in a fixed neighborhood of $\mathbf t=0$,
\[
 \log\E e^{\sum_i t_iM_{i,\varepsilon}}
 =\frac{P_{\mathbf L}(e^{t_1},\ldots,e^{t_r})-1}{\varepsilon}+O_{\mathbf L}(1)
\]
with locally bounded derivatives of every fixed order.  For any nonzero
multi-index $\boldsymbol\alpha$, differentiation at the origin gives
\[
 \frac1{|G|}\sum_{g\in G}
 \prod_{i\in S(\boldsymbol\alpha)}m_i(g)
 =\frac{|\operatorname{Gal}(L/L_{S(\boldsymbol\alpha)})|}{|G|}
 =\frac1{[L_{S(\boldsymbol\alpha)}:\mathbb Q]},
\]
which proves \eqref{eq:Galois-cumulant-compositum} and
\eqref{eq:Galois-means-covariance}.  Since $L_i$ and $L_j$ are Galois over
$\mathbb Q$,
\[
 [L_iL_j:\mathbb Q]\,[L_i\cap L_j:\mathbb Q]
 =[L_i:\mathbb Q][L_j:\mathbb Q],
\]
and \eqref{eq:Galois-intersection-recovery} follows.  The same substitution of
\eqref{eq:Galois-cumulant-compositum} yields
\eqref{eq:Galois-entanglement-factor}.
\end{proof}

For the moving-profile analysis, write
\[
 \mathcal V_{\mathbf L}
 :=\{\mathbf v\in\{0,1\}^r\setminus\{\mathbf0\}:\pi_{\mathbf v}>0\},
 \qquad
 \Lambda_{\mathbf L}:=\operatorname{span}_{\mathbb Z}\mathcal V_{\mathbf L},
\]
and let
\[
 \mathscr C_{\mathbf L}
 :=\operatorname{cone}(\mathcal V_{\mathbf L})\subset\mathbb R^r.
\]

\begin{lemma}
\label{lem:Galois-semigroup-saturation}
Let
\[
 \Gamma_{\mathbf L}:=\mathbb Z_{\ge0}\mathcal V_{\mathbf L},\qquad
 \overline\Gamma_{\mathbf L}:=\mathscr C_{\mathbf L}\cap\Lambda_{\mathbf L}.
\]
There exists a conductor vector $\mathbf c_{\mathbf L}\in\Gamma_{\mathbf L}$
such that
\[
 \mathbf c_{\mathbf L}+\overline\Gamma_{\mathbf L}
 \subseteq\Gamma_{\mathbf L}.
\]
Consequently, for every compact $K\Subset\mathscr C_{\mathbf L}^{\circ}$
there is $\varepsilon_0(K)>0$ such that, whenever
$0<\varepsilon<\varepsilon_0(K)$,
$\mathbf q\in\Lambda_{\mathbf L}\cap\mathbb Z_{\ge0}^r$, and
$\varepsilon\mathbf q\in K$, one has $\mathbf q\in\Gamma_{\mathbf L}$.
\end{lemma}

\begin{proof}
Every $\mathbf y\in\overline\Gamma_{\mathbf L}$ admits a
representation $\mathbf y=\sum_{\mathbf v\in\mathcal V_{\mathbf L}}
a_{\mathbf v}\mathbf v$ with $a_{\mathbf v}\ge0$.  Subtracting the integer
parts of the coefficients shows that
\[
 \mathbf y=\boldsymbol\gamma+\mathbf b,
 \qquad \boldsymbol\gamma\in\Gamma_{\mathbf L},
\]
where $\mathbf b$ belongs to the bounded zonotope
$\{\sum_{\mathbf v}t_{\mathbf v}\mathbf v:0\le t_{\mathbf v}<1\}$ and to
$\Lambda_{\mathbf L}$.  Hence only finitely many such remainders
$\mathbf b$ occur.  Since each lies in the group generated by
$\mathcal V_{\mathbf L}$, write
$\mathbf b=\boldsymbol\alpha_{\mathbf b}-\boldsymbol\delta_{\mathbf b}$
with $\boldsymbol\alpha_{\mathbf b},\boldsymbol\delta_{\mathbf b}
\in\Gamma_{\mathbf L}$.  Taking
$\mathbf c_{\mathbf L}:=\sum_{\mathbf b}\boldsymbol\delta_{\mathbf b}$ gives
$\mathbf c_{\mathbf L}+\mathbf b\in\Gamma_{\mathbf L}$ for every remainder,
and therefore
$\mathbf c_{\mathbf L}+\overline\Gamma_{\mathbf L}\subseteq
\Gamma_{\mathbf L}$.

For the final assertion, compactness and
$K\Subset\mathscr C_{\mathbf L}^{\circ}$ imply
$\mathbf x-\varepsilon\mathbf c_{\mathbf L}\in\mathscr C_{\mathbf L}$
uniformly for $\mathbf x\in K$ and all sufficiently small $\varepsilon$.
Thus, if $\mathbf x=\varepsilon\mathbf q$ with
$\mathbf q\in\Lambda_{\mathbf L}$, then
$\mathbf q-\mathbf c_{\mathbf L}\in
\mathscr C_{\mathbf L}\cap\Lambda_{\mathbf L}
=\overline\Gamma_{\mathbf L}$, and the conductor inclusion yields
$\mathbf q\in\Gamma_{\mathbf L}$.
\end{proof}

\paragraph{Two basic full-rank situations.}
The rank hypothesis in the moving-profile theorem is automatic in several
standard arithmetic configurations.  If the nontrivial Galois extensions
$L_1,\ldots,L_r$ are linearly disjoint over $\mathbb Q$, then
$\operatorname{Gal}(L_1\cdots L_r/\mathbb Q)\simeq
\prod_i\operatorname{Gal}(L_i/\mathbb Q)$.  Choosing the identity in the
$i$th factor and nonidentity elements in the others realizes the jump vector
$e_i$, so $\Lambda_{\mathbf L}$ has rank $r$.  For a proper tower
$L_1\subsetneq L_2$ of finite Galois extensions, the group-level patterns
$(1,0)$ and $(1,1)$ both occur, and the associated two-dimensional marking
lattice again has full rank.  Degeneracy can occur, for example, when two
marked fields coincide; then the corresponding coordinates are identical and
the rank drops.

\begin{theorem}
\label{thm:Galois-moving-profile}
Assume that $\Lambda_{\mathbf L}$ has rank $r$, and put
\[
 h_{\mathbf L}:=[\mathbb Z^r:\Lambda_{\mathbf L}].
\]
For $\mathbf t\in\mathbb R^r$ define
\begin{equation}\label{eq:Galois-cumulant-potential}
 Q_{\mathbf L}(\mathbf t)
 :=P_{\mathbf L}(e^{t_1},\ldots,e^{t_r})-1
 =\sum_{\mathbf v\in\mathcal V_{\mathbf L}}
   \pi_{\mathbf v}\bigl(e^{\langle\mathbf v,\mathbf t\rangle}-1\bigr).
\end{equation}
For every $\mathbf x\in\mathscr C_{\mathbf L}^{\circ}$ there is a unique
$\mathbf t(\mathbf x)\in\mathbb R^r$ satisfying
\begin{equation}\label{eq:Galois-saddle-system}
 \nabla Q_{\mathbf L}(\mathbf t(\mathbf x))=\mathbf x.
\end{equation}
Set
\[
 \boldsymbol\zeta(\mathbf x):=e^{\mathbf t(\mathbf x)},
 \qquad
 \Sigma_{\mathbf L}(\mathbf x)
 :=\nabla^2Q_{\mathbf L}(\mathbf t(\mathbf x))
 =\sum_{\mathbf v\in\mathcal V_{\mathbf L}}
   \pi_{\mathbf v}\boldsymbol\zeta(\mathbf x)^{\mathbf v}
   \mathbf v\mathbf v^{\!T},
\]
and
\begin{equation}\label{eq:Galois-rate-function}
 I_{\mathbf L}(\mathbf x)
 :=\langle\mathbf x,\mathbf t(\mathbf x)\rangle
   -Q_{\mathbf L}(\mathbf t(\mathbf x)).
\end{equation}
Then $\Sigma_{\mathbf L}(\mathbf x)$ is positive definite and
$I_{\mathbf L}$ is the strictly convex Legendre transform of
$Q_{\mathbf L}$ on $\mathscr C_{\mathbf L}^{\circ}$.

Write
\[
 F_{\mathbf L,\mathbf q}(\sigma)
 :=[\mathbf z^{\mathbf q}]\mathcal F_{\mathbf L}(\sigma,\mathbf z).
\]
For every compact $K\Subset\mathscr C_{\mathbf L}^{\circ}$, uniformly for
$\mathbf q\in\mathbb Z_{\ge0}^r\cap\Lambda_{\mathbf L}$ with
$\mathbf x:=\varepsilon\mathbf q\in K$ (and hence, by
Lemma~\ref{lem:Galois-semigroup-saturation}, $\mathbf q\in\Gamma_{\mathbf L}$
for all sufficiently small $\varepsilon$),
\begin{align}
 F_{\mathbf L,\mathbf q}(1+\varepsilon)
 &=\frac{h_{\mathbf L}\,\varepsilon^{r/2}}
 {(2\pi)^{r/2}\sqrt{\det\Sigma_{\mathbf L}(\mathbf x)}}
 \mathcal C_{\mathbf L,0}(\boldsymbol\zeta(\mathbf x))
 \boldsymbol\zeta(\mathbf x)^{-\mathbf q}\notag\\
 &\qquad\times
 \exp\!\left(\frac{P_{\mathbf L}(\boldsymbol\zeta(\mathbf x))}{\varepsilon}\right)
 \bigl(1+O_{\mathbf L,K}(\varepsilon)\bigr).
 \label{eq:Galois-local-coefficient}
\end{align}
If $\mathbf q\notin\Lambda_{\mathbf L}$, then
$F_{\mathbf L,\mathbf q}(1+\varepsilon)=0$ identically.

Consequently, under the probability measure
\eqref{eq:Galois-weighted-ensemble}, uniformly in the same moving regime,
\begin{align}
 \Prob(\mathbf M_\varepsilon=\mathbf q)
 &=\frac{h_{\mathbf L}\,\varepsilon^{r/2}}
 {(2\pi)^{r/2}\sqrt{\det\Sigma_{\mathbf L}(\mathbf x)}}
 \frac{\mathcal C_{\mathbf L,0}(\boldsymbol\zeta(\mathbf x))}
      {\mathcal C_{\mathbf L,0}(\mathbf1)}
 \exp\!\left(-\frac{I_{\mathbf L}(\mathbf x)}{\varepsilon}\right)
 \bigl(1+O_{\mathbf L,K}(\varepsilon)\bigr).
 \label{eq:Galois-local-probability}
\end{align}
In particular,
\begin{equation}\label{eq:Galois-local-large-deviation}
 \varepsilon\log\Prob(\mathbf M_\varepsilon=\mathbf q)
 =-I_{\mathbf L}(\varepsilon\mathbf q)+O_{\mathbf L,K}
   (\varepsilon|\log\varepsilon|).
\end{equation}
Thus the leading local moving-support profile, including its large-deviation
geometry, is determined by the same polynomial $P_{\mathbf L}$ whose mixed
coefficients are the reciprocal compositum degrees.
\end{theorem}

\begin{proof}
Because the nonzero jump vectors span $\mathbb R^r$, the Hessian
\[
 \nabla^2Q_{\mathbf L}(\mathbf t)
 =\sum_{\mathbf v\in\mathcal V_{\mathbf L}}
   \pi_{\mathbf v}e^{\langle\mathbf v,\mathbf t\rangle}
   \mathbf v\mathbf v^{\!T}
\]
is positive definite for every $\mathbf t$.  Hence $Q_{\mathbf L}$ is strictly
convex.  If $\mathbf x\in\mathscr C_{\mathbf L}^{\circ}$, the function
$Q_{\mathbf L}(\mathbf t)-\langle\mathbf x,\mathbf t\rangle$ is coercive.
Indeed, if $\|\mathbf t_n\|\to\infty$, pass to a subsequence for which
$\mathbf t_n/\|\mathbf t_n\|\to\mathbf d$.  If
$\langle\mathbf v,\mathbf d\rangle>0$ for some
$\mathbf v\in\mathcal V_{\mathbf L}$, the corresponding exponential term
forces the expression to $+\infty$.  Otherwise $\mathbf d$ lies in the polar
cone of $\mathscr C_{\mathbf L}$; since $\mathbf x$ is an interior point and
$\mathbf d\ne0$, one has $\langle\mathbf x,\mathbf d\rangle<0$, so the
linear term $-\langle\mathbf x,\mathbf t_n\rangle$ tends to $+\infty$
along the subsequence.  It therefore has a unique minimizer, characterized by
\eqref{eq:Galois-saddle-system}.  This also proves the stated
Legendre-duality assertions.

Fix $K\Subset\mathscr C_{\mathbf L}^{\circ}$.  The vectors
$\mathbf t(\mathbf x)$ and $\boldsymbol\zeta(\mathbf x)$ remain in compact
sets as $\mathbf x$ ranges over $K$.  Cauchy's formula on the torus
$z_j=\zeta_j e^{i\theta_j}$ gives
\begin{align}
 F_{\mathbf L,\mathbf q}(1+\varepsilon)
 =\frac{\boldsymbol\zeta^{-\mathbf q}}{(2\pi)^r}
 \int_{[-\pi,\pi]^r}
  \mathcal F_{\mathbf L}(1+\varepsilon,
       \boldsymbol\zeta e^{i\boldsymbol\theta})
  e^{-i\langle\mathbf q,\boldsymbol\theta\rangle}\,
  \dd\boldsymbol\theta.
 \label{eq:Galois-Cauchy-saddle}
\end{align}
By Theorem~\ref{thm:Galois-entanglement-support}, uniformly on all tori that
occur here,
\[
 \mathcal F_{\mathbf L}(1+\varepsilon,\mathbf z)
 =e^{P_{\mathbf L}(\mathbf z)/\varepsilon}
  \bigl(\mathcal C_{\mathbf L,0}(\mathbf z)+O_{\mathbf L,K}(\varepsilon)\bigr).
\]
The limiting factor is positive on $(0,\infty)^r$, as is clear from the
canonical-product representation in the proof of
Theorem~\ref{thm:Galois-entanglement-support}.

Let
\[
 \Theta_{\mathbf L}
 :=\{\boldsymbol\theta\in(\mathbb R/2\pi\mathbb Z)^r:
      \langle\mathbf v,\boldsymbol\theta\rangle\in2\pi\mathbb Z
      \text{ for all }\mathbf v\in\Lambda_{\mathbf L}\}.
\]
This finite dual group has cardinality $h_{\mathbf L}$.  Since
$P_{\mathbf L}(\mathbf z)=\sum_{\mathbf v}\pi_{\mathbf v}\mathbf z^{\mathbf v}$,
\[
 \Re P_{\mathbf L}(\boldsymbol\zeta e^{i\boldsymbol\theta})
 \le P_{\mathbf L}(\boldsymbol\zeta),
\]
with equality exactly for $\boldsymbol\theta\in\Theta_{\mathbf L}$.  Choose
fixed disjoint neighborhoods of these finitely many points.  Since
$\boldsymbol\zeta(\mathbf x)$ ranges over a compact subset of
$(0,\infty)^r$, continuity on the compact product of this saddle set with the
complement of those neighborhoods gives a uniform gap
\[
 P_{\mathbf L}(\boldsymbol\zeta(\mathbf x))
 -\Re P_{\mathbf L}(\boldsymbol\zeta(\mathbf x)e^{i\boldsymbol\theta})
 \ge \delta_K>0.
\]
Thus the corresponding part of \eqref{eq:Galois-Cauchy-saddle} is uniformly
exponentially smaller than the main term.

Near $\boldsymbol\theta_0\in\Theta_{\mathbf L}$ write
$\boldsymbol\theta=\boldsymbol\theta_0+\mathbf u$.  Since
$\varepsilon\mathbf q=\mathbf x=\nabla Q_{\mathbf L}(\mathbf t(\mathbf x))$,
Taylor expansion yields
\begin{align*}
 &P_{\mathbf L}(\boldsymbol\zeta e^{i\mathbf u})
   -P_{\mathbf L}(\boldsymbol\zeta)
   -i\langle\mathbf x,\mathbf u\rangle\\
 &\qquad=-\frac12\mathbf u^T\Sigma_{\mathbf L}(\mathbf x)\mathbf u
   +O_{\mathbf L,K}(\|\mathbf u\|^3).
\end{align*}
Because $K$ is compact, the saddle set $\{\mathbf t(\mathbf x):\mathbf x\in K\}$
is compact; hence all third and fourth derivatives of the phase and the
required derivatives of the canonical amplitude are uniformly bounded on a
fixed complex neighborhood of the relevant tori.  Together with the uniform
minor-arc gap above and the positive lower bound for the smallest eigenvalue
of $\Sigma_{\mathbf L}(\mathbf x)$ on $K$, this supplies a uniform Gaussian
majorant.  The standard uniform multivariate Laplace expansion (see, for
example, Pemantle--Wilson \cite[Ch.~5]{PemantleWilson2013}) therefore gives
\[
 \int
 e^{[P_{\mathbf L}(\boldsymbol\zeta e^{i\mathbf u})
      -P_{\mathbf L}(\boldsymbol\zeta)
      -i\langle\mathbf x,\mathbf u\rangle]/\varepsilon}
 \mathcal C_{\mathbf L,0}(\boldsymbol\zeta e^{i\mathbf u})\,
 \dd\mathbf u
 =\frac{(2\pi\varepsilon)^{r/2}}
       {\sqrt{\det\Sigma_{\mathbf L}(\mathbf x)}}
  \mathcal C_{\mathbf L,0}(\boldsymbol\zeta)
  \bigl(1+O_{\mathbf L,K}(\varepsilon)\bigr).
\]
The $O(\varepsilon)$ remainder is uniform on $K$; the potential
$O(\sqrt\varepsilon)$ terms are odd after the rescaling
$\mathbf u=\sqrt\varepsilon\,\mathbf y$ and integrate to zero in the first
Laplace correction.  The uniform $O_{\mathbf L,K}(\varepsilon)$ amplitude error in the
preceding approximation to $\mathcal F_{\mathbf L}$ is handled by the same
Gaussian majorant.  Moreover,
$\mathcal C_{\mathbf L,0}(\boldsymbol\zeta(\mathbf x))$ is positive and,
by compactness of $K$, bounded uniformly away from zero.  The integrated
amplitude error is therefore $O_{\mathbf L,K}(\varepsilon)$ times the Gaussian
main scale and hence contributes a relative
$O_{\mathbf L,K}(\varepsilon)$ term.

Both $P_{\mathbf L}$ and $\mathcal C_{\mathbf L,0}$ are invariant under the
dual group $\Theta_{\mathbf L}$.  Indeed, for every $\varepsilon>0$ the Euler
product $\mathcal F_{\mathbf L}(1+\varepsilon,\mathbf z)$ depends on
$\mathbf z$ only through monomials with exponents in $\Lambda_{\mathbf L}$,
and $P_{\mathbf L}$ has the same invariance.  Hence the normalized entire
function in \eqref{eq:Galois-canonical-limit} is $\Theta_{\mathbf L}$-invariant
for every $\varepsilon$, and local uniform convergence passes this invariance
to $\mathcal C_{\mathbf L,0}$.  If
$\mathbf q\in\Lambda_{\mathbf L}$, then
$e^{-i\langle\mathbf q,\boldsymbol\theta_0\rangle}=1$ for every
$\boldsymbol\theta_0\in\Theta_{\mathbf L}$, so all $h_{\mathbf L}$ saddle
contributions are equal.  This proves \eqref{eq:Galois-local-coefficient}.
If $\mathbf q\notin\Lambda_{\mathbf L}$, character orthogonality on the dual
group shows that the coefficient is zero; equivalently, every monomial in the
Euler product has exponent in $\Lambda_{\mathbf L}$.

Finally,
\[
 \mathcal F_{\mathbf L}(1+\varepsilon,\mathbf1)
 =e^{1/\varepsilon}\mathcal C_{\mathbf L,0}(\mathbf1)
  \bigl(1+O_{\mathbf L}(\varepsilon)\bigr),
\]
and
\[
 \boldsymbol\zeta^{-\mathbf q}
 \exp\!\left(\frac{P_{\mathbf L}(\boldsymbol\zeta)-1}{\varepsilon}\right)
 =\exp\!\left(-\frac{I_{\mathbf L}(\mathbf x)}{\varepsilon}\right).
\]
Division proves \eqref{eq:Galois-local-probability}; taking logarithms gives
\eqref{eq:Galois-local-large-deviation}.
\end{proof}

\begin{corollary}\label{cor:Galois-local-CLT}
Under the hypotheses of Theorem~\ref{thm:Galois-moving-profile}, let
\[
 \bar{\mathbf x}:=\nabla Q_{\mathbf L}(\mathbf0)
 =\left(\frac1{[L_1:\mathbb Q]},\ldots,
        \frac1{[L_r:\mathbb Q]}\right)
\]
and
\[
 \Sigma_{\mathbf L,0}:=\nabla^2Q_{\mathbf L}(\mathbf0),
 \qquad
 (\Sigma_{\mathbf L,0})_{ij}
 =\frac1{[L_iL_j:\mathbb Q]},
\]
where $L_iL_i=L_i$.  If
$\mathbf q\in\mathbb Z_{\ge0}^r\cap\Lambda_{\mathbf L}$ and
\[
 \mathbf y_\varepsilon
 :=\sqrt\varepsilon\left(\mathbf q-\frac{\bar{\mathbf x}}{\varepsilon}\right)
\]
remains in a fixed compact subset of $\mathbb R^r$, then
\begin{align}
 \Prob(\mathbf M_\varepsilon=\mathbf q)
 &=\frac{h_{\mathbf L}\,\varepsilon^{r/2}}
 {(2\pi)^{r/2}\sqrt{\det\Sigma_{\mathbf L,0}}}
 \exp\!\left(-\frac12
  \mathbf y_\varepsilon^T\Sigma_{\mathbf L,0}^{-1}
  \mathbf y_\varepsilon\right)
 \bigl(1+O_{\mathbf L}(\sqrt\varepsilon)\bigr).
 \label{eq:Galois-local-CLT}
\end{align}
Thus the covariance matrix in the Gaussian window has entries
$[L_iL_j:\mathbb Q]^{-1}/\varepsilon$, in agreement with the mixed-cumulant
formula \eqref{eq:Galois-cumulant-compositum}.
\end{corollary}

\begin{proof}
At $\mathbf x=\bar{\mathbf x}$ the saddle is $\mathbf t=\mathbf0$ and
$I_{\mathbf L}(\bar{\mathbf x})=0$.  Taylor expansion of the Legendre transform
gives
\[
 I_{\mathbf L}(\bar{\mathbf x}+\mathbf h)
 =\frac12\mathbf h^T\Sigma_{\mathbf L,0}^{-1}\mathbf h
  +O_{\mathbf L}(\|\mathbf h\|^3).
\]
Apply Theorem~\ref{thm:Galois-moving-profile} with
$\mathbf h=\sqrt\varepsilon\,\mathbf y_\varepsilon$ and expand the smooth
prefactor at $\bar{\mathbf x}$.
\end{proof}

\begin{example}
\label{ex:three-quadratic-entanglement}
Take
\[
 L_1=\mathbb Q(\sqrt2),\qquad
 L_2=\mathbb Q(\sqrt3),\qquad
 L_3=\mathbb Q(\sqrt6).
\]
Every pair has trivial intersection, so each normalized pairwise covariance in
\eqref{eq:Galois-intersection-recovery} tends to $1$.  However
$L_1L_2L_3=\mathbb Q(\sqrt2,\sqrt3)$ has degree $4$, whereas
$\prod_i[L_i:\mathbb Q]=8$.  Hence the normalized third mixed cumulant in
\eqref{eq:Galois-entanglement-factor} tends to $2$.  Thus the higher joint
support statistics detect a Galois relation invisible to all pairwise
intersection degrees.
\end{example}

\section{Bernoulli--Barnes kernels and normalized remainders}
\label{sec:centered}

This section develops the explicit special-function model used in the geometric
realization.  We first represent the centered Bernoulli--Barnes kernel as the
Laplace transform of a positive random variable, and then use that representation
to construct normalized remainders with a finite family of Mellin poles.  These
two steps supply the analytic data needed for the Barnes--Epstein factorization
and the pole-collision results below.

\subsection{A positive Bernoulli--Barnes representation}

Let
\[
  \ba=(a_1,\ldots,a_d)\in(\R\setminus\{0\})^d,
  \qquad A=\sum_{j=1}^d a_j.
\]
For the Bernoulli--Barnes polynomials (see, e.g., \cite{BayadBeck2014,BayadKim2013})
\[
 \frac{z^d e^{xz}}{\prod_{j=1}^d(e^{a_jz}-1)}
 =\sum_{n\ge0}\BB_n(x;\ba)\frac{z^n}{n!},
\]
define the normalized centered coefficients
\begin{equation}\label{eq:Cn-definition}
 \Ca_n(\ba)
 :=\left(\prod_{j=1}^d a_j\right)\BB_n(A/2;\ba),
 \qquad
 \Phi_{\ba}(z):=\sum_{n\ge0}\Ca_n(\ba)\frac{z^n}{n!}.
\end{equation}
At the center one has
\begin{equation}\label{eq:centered-product}
 \Phi_{\ba}(z)
 =\prod_{j=1}^d\frac{a_jz}{2\sinh(a_jz/2)}.
\end{equation}
Thus \(\Ca_{2n+1}(\ba)=0\), and \(\Ca_n(\ba)=\Ca_n(|\ba|)\).

For \(y\ge0\) put
\begin{equation}\label{eq:F-definition}
 \Fa_{\ba}(y):=\Phi_{\ba}(\sqrt y).
\end{equation}
Let \(\{E_{j,k}\}\) be independent standard exponential random variables and set
\begin{equation}\label{eq:X-definition}
 \lambda_{j,k}:=\frac{a_j^2}{4\pi^2k^2},\qquad
 X_{\ba}:=\sum_{j=1}^d\sum_{k=1}^{\infty}\lambda_{j,k}E_{j,k},
 \qquad M_n(\ba):=\E[X_{\ba}^n].
\end{equation}

The following theorem provides the positive probabilistic representation on
which the remainder and Mellin constructions of the subsequent sections are
based.

\begin{theorem}
\label{thm:master-probability}
The series in \eqref{eq:X-definition} converges almost surely and in \(L^1\),
all positive moments are finite, and
\begin{equation}\label{eq:Laplace-main}
 \Fa_{\ba}(y)
 =\prod_{j=1}^d\prod_{k=1}^{\infty}(1+\lambda_{j,k}y)^{-1}
 =\E[e^{-yX_{\ba}}].
\end{equation}
Moreover,
\begin{equation}\label{eq:C-moment}
 \frac{\Ca_{2n}(\ba)}{(2n)!}
 =\frac{(-1)^n}{n!}M_n(\ba),
 \qquad (-1)^n\Ca_{2n}(\ba)>0\quad(n\ge1).
\end{equation}
In particular, $X_{\ba}$ is a generalized gamma convolution with discrete
Thorin measure
\begin{equation}\label{eq:Thorin-measure-final}
 \mathcal U_{\ba}
 =\sum_{j=1}^d\sum_{k\ge1}
  \delta_{4\pi^2k^2/a_j^2};
\end{equation}
see \cite{Bondesson1992} for the general GGC framework.
\end{theorem}

\begin{proof}
Euler's product
\[
 \frac{\sinh w}{w}
 =\prod_{k\ge1}\left(1+\frac{w^2}{\pi^2k^2}\right)
\]
with \(w=a_j\sqrt y/2\) gives the product in \eqref{eq:Laplace-main}.
Since
\[
 \sum_{j,k}\lambda_{j,k}=\frac1{24}\sum_j a_j^2<\infty,
\]
the nonnegative series defining $X_{\ba}$ converges almost surely and in
$L^1$.  If $0<t<\lambda_{\max}^{-1}$, where
$\lambda_{\max}=\max_{j,k}\lambda_{j,k}$, then
$\prod_{j,k}(1-t\lambda_{j,k})^{-1}<\infty$; hence the moment generating
function is finite near the origin and all positive moments are finite.
Independence and \(\E(e^{-y\lambda E})=(1+\lambda y)^{-1}\) give the Laplace
identity.  Since
\[
 \log \Fa_{\ba}(y)
 =-\sum_{j,k}\log\!\left(1+\frac{y}{4\pi^2k^2/a_j^2}\right)
\]
and $\sum_{j,k}(1+4\pi^2k^2/a_j^2)^{-1}<\infty$, this is the Thorin
representation with measure \eqref{eq:Thorin-measure-final}.  Expanding the
Laplace identity at the origin yields \eqref{eq:C-moment}.
\end{proof}

\subsection{Normalized remainders and Mellin poles}
\label{sec:remainders}

The normalized remainders introduced in this subsection provide the analytic
source of the additional Mellin poles responsible for the transition proved
later.  The positive representation from the preceding subsection turns these
remainders into completely monotone kernels and yields a finite inverse-power
expansion whose poles subsequently collide with the Epstein pole.  This is the
mechanism that raises the critical pole multiplicity in the resonant dimensions.

For background on the normalized-remainder terminology and related
developments, see \cite{QiSurvey2025}.  For \(m\ge0\) define
\begin{equation}\label{eq:Remainder-def}
 \Ra_{m,\ba}(y)
 :=\Fa_{\ba}(y)-\sum_{n=0}^{m}\Ca_{2n}(\ba)\frac{y^n}{(2n)!},
\end{equation}
and
\begin{equation}\label{eq:T-def}
 \Ta_{m,\ba}(y)
 :=\frac{(2m+2)!}{\Ca_{2m+2}(\ba)y^{m+1}}\Ra_{m,\ba}(y),
 \qquad \Ta_{m,\ba}(0):=1.
\end{equation}
For a nonnegative random variable \(X\), let \(X^{[r]}\) denote the \(r\)th
power-biased law,
\begin{equation}\label{eq:size-biased}
 \E h(X^{[r]})=\frac{\E[X^rh(X)]}{\E X^r}.
\end{equation}

\begin{lemma}\label{lem:rapid-decay}
For every \(r,L\in\mathbb Z_{\ge0}\) there is \(C_{r,L}(\ba)>0\) such that
\begin{equation}\label{eq:rapid-decay}
 \E[X_{\ba}^r e^{-yX_{\ba}}]
 \le C_{r,L}(\ba)(1+y)^{-L},\qquad y\ge0.
\end{equation}
\end{lemma}

\begin{proof}
For \(y\ge1\), \(x^re^{-yx}\le C_r y^{-r}e^{-yx/2}\), so
\[
 \E[X_{\ba}^re^{-yX_{\ba}}]\le C_r y^{-r}\Fa_{\ba}(y/2).
\]
If \(L>r\), retain \(L-r\) distinct factors in \eqref{eq:Laplace-main}; all
remaining factors are at most one, giving \(\Fa_{\ba}(y/2)=O(y^{-(L-r)})\).
If \(L\le r\), use \(\Fa_{\ba}(y/2)\le1\).  The interval \(0\le y\le1\)
is immediate.
\end{proof}

The next theorem gives both the positive representation and the precise
finite asymptotic expansion that will determine the resonance set.

\begin{theorem}\label{thm:normalized-remainder}
Let \(U_m\sim\operatorname{Beta}(1,m+1)\) be independent of
\(X_{\ba}^{[m+1]}\), and put \(Z_{m,\ba}=U_mX_{\ba}^{[m+1]}\).  Then
\begin{equation}\label{eq:T-Laplace}
 \Ta_{m,\ba}(y)=\E[e^{-yZ_{m,\ba}}]
 =(m+1)\int_0^1(1-u)^m
 \frac{\E[X_{\ba}^{m+1}e^{-uyX_{\ba}}]}{M_{m+1}(\ba)}\dd u.
\end{equation}
In particular \(\Ta_{m,\ba}\) is completely monotone.  For every fixed
\(m\ge0\) and every integer \(L>m+1\),
\begin{align}
 \Ta_{m,\ba}(y)
 &=\sum_{k=1}^{m+1}c_{m,k}(\ba)y^{-k}+O_L(y^{-L}),
 \label{eq:T-infty-coeff}\\
 c_{m,k}(\ba)
 &=(-1)^{k-1}\frac{(m+1)!}{(m-k+1)!}
 \frac{M_{m-k+1}(\ba)}{M_{m+1}(\ba)}.
 \label{eq:T-infty-cmk}
\end{align}
\end{theorem}

\begin{proof}
For $x,y\ge0$, Taylor's formula with integral remainder gives
\[
 e^{-yx}-\sum_{n=0}^{m}\frac{(-yx)^n}{n!}
 =\frac{(-yx)^{m+1}}{m!}
   \int_0^1(1-u)^m e^{-uyx}\dd u.
\]
Taking $x=X_{\ba}$ and then expectation is justified by positivity after
multiplication by $(-1)^{m+1}$, or directly by the finiteness of all moments.
Using \eqref{eq:C-moment} to identify the Taylor polynomial of
$\Fa_{\ba}(y)$, we obtain
\[
 \Ra_{m,\ba}(y)
 =\frac{(-1)^{m+1}y^{m+1}}{m!}
  \int_0^1(1-u)^m
  \E[X_{\ba}^{m+1}e^{-uyX_{\ba}}]\dd u.
\]
Since
$\Ca_{2m+2}/(2m+2)!=(-1)^{m+1}M_{m+1}/(m+1)!$, the normalization
\eqref{eq:T-def} reduces this identity exactly to
\[
 \Ta_{m,\ba}(y)
 =(m+1)\int_0^1(1-u)^m
 \frac{\E[X_{\ba}^{m+1}e^{-uyX_{\ba}}]}{M_{m+1}}\dd u.
\]
Now $U_m$ has density $(m+1)(1-u)^m\mathbf 1_{(0,1)}(u)$, while by the
definition of power bias,
\[
 \E[e^{-yU_mX_{\ba}^{[m+1]}}\mid U_m=u]
 =\frac{\E[X_{\ba}^{m+1}e^{-uyX_{\ba}}]}{M_{m+1}}.
\]
Averaging over $U_m$ proves the Laplace representation and hence complete
monotonicity.

For the large-$y$ expansion put
$f_m(v)=\E[X_{\ba}^{m+1}e^{-vX_{\ba}}]$ and substitute $v=uy$:
\begin{equation}\label{eq:T-expansion-proof-detail}
 \Ta_{m,\ba}(y)
 =\frac{m+1}{yM_{m+1}}
 \int_0^y\left(1-\frac vy\right)^m f_m(v)\dd v.
\end{equation}
Lemma~\ref{lem:rapid-decay} implies, for every prescribed $L$, that replacing
the upper limit $y$ by $\infty$ in \eqref{eq:T-expansion-proof-detail}
changes the result by $O_L(y^{-L})$.  Because the factor in parentheses is a
polynomial of degree $m$,
\[
 \left(1-\frac vy\right)^m
 =\sum_{\ell=0}^{m}(-1)^\ell\binom m\ell v^\ell y^{-\ell}.
\]
Tonelli's theorem and
$\int_0^\infty v^\ell e^{-vX}\dd v=\ell!X^{-\ell-1}$ give
\[
 \int_0^\infty v^\ell f_m(v)\dd v
 =\ell!\E[X_{\ba}^{m-\ell}]
 =\ell!M_{m-\ell},\qquad 0\le\ell\le m.
\]
Substitution into \eqref{eq:T-expansion-proof-detail}, followed by the change
of index $k=\ell+1$, yields precisely
\eqref{eq:T-infty-coeff} and \eqref{eq:T-infty-cmk}.
\end{proof}

\section{Barnes--Epstein realization and exact-order torsion}
\label{sec:arithmetic-mellin}

We now connect the abstract Mellin-resonance theory to elliptic geometry.
The Bernoulli--Barnes moment transform factors through a Barnes zeta
function, while lattice summation introduces an Epstein factor.  This
separation makes the critical pole explicit.  Finite Fourier orthogonality
then identifies exact-order torsion averaging with the operator
$\mathscr P_{N,2}$, so the local arithmetic filtration of
Section~\ref{sec:abstract-resonance} becomes a geometric statement about
elliptic torsion.

\paragraph{Convention on Mellin transforms.}
Whenever the defining integral converges, $\widehat H(s)$ denotes the ordinary
Mellin transform.  The same notation is used for its meromorphic continuation
when the continuation is under discussion.  Unless ordinary integrability is
explicitly asserted, a critical value or derivative refers to this continued
Mellin transform.  The final subsection below gives a strip-margin criterion
and verifies when the critical derivatives in the exact-order families are
genuine logarithmic Mellin moments.

\paragraph{Regularization convention.}
Three analytically different objects occur below and will be kept distinct.
The theta kernels are defined by ordinary absolutely convergent lattice sums.
The two-dimensional normalized-remainder kernel is instead introduced first by
the Mellin-regularized inverse transform \eqref{eq:twisted-lattice-remainder},
because the unsubtracted lattice sum is not absolutely convergent.  In higher
dimension, after subtracting the finitely many inverse-power terms required by
\eqref{eq:subtracted-Mellin-v15}, the renormalized lattice kernel is again an
ordinary absolutely convergent sum on a nonempty strip.  Whenever exact-order
cancellation improves the endpoint behavior enough to recover an ordinary
Mellin integral, this is stated explicitly rather than inferred from
meromorphic continuation alone.

Throughout this section put
\[
 \bb=(|a_1|,\ldots,|a_d|),\qquad
 B=\sum_{j=1}^d|a_j|,\qquad D=\prod_{j=1}^d|a_j|,
\]
and write
\[
 \Fa_{\ba}(y)=\E(e^{-y\Xa}),\qquad
 \mu_{\ba}(q)=\E(\Xa^q).
\]
We initially use $\mu_{\ba}(q)=\E(X_{\ba}^q)$ whenever the expectation is finite.  The results below show that $z\mapsto\mu_{\ba}(z)$ is in fact entire, which is the form needed for the Mellin arguments.

\subsection{The Bernoulli--Barnes/Barnes--Epstein bridge}

For $\bb=(b_1,\ldots,b_d)\in(0,\infty)^d$, $x>0$, and initially
$\Re s>d$, write
\[
 \zeta_d(s,x;\bb)
 :=\sum_{\bm m\in\Z_{\ge0}^d}
 (x+m_1b_1+\cdots+m_db_d)^{-s}
\]
for the Barnes zeta function, continued meromorphically in the usual way; see
\cite{Barnes1904,BayadBeck2014}.

\begin{theorem}
\label{thm:Barnes-Epstein-bridge}
For $\Re s>0$,
\begin{equation}\label{eq:negative-moment}
 \mu_{\ba}(-s)
 =2D\,\frac{\Gamma(2s+d)}{\Gamma(s)}
 \zeta_d\!\left(2s+d,\frac B2;\bb\right).
\end{equation}
Let $\tau=u+iv\in\mathbb H$ and
\[
 Q_\tau(m,n)=\frac{|m+n\tau|^2}{v},\qquad
 \Lambda_\tau(s)=\pi^{-s}\Gamma(s)
 \sum_{(m,n)\ne(0,0)}Q_\tau(m,n)^{-s}.
\]
The Epstein zeta function $\Lambda_\tau$ was introduced in \cite{Epstein1903}.
Define
\begin{equation}\label{eq:theta-GGC}
 \Theta_{\ba}(\tau;t)
 :=\sum_{(m,n)\ne(0,0)}
 \Fa_{\ba}\!\left(\pi tQ_\tau(m,n)\right).
\end{equation}
Then, initially for $\Re s>1$,
\begin{equation}\label{eq:Barnes-Epstein}
 \int_0^\infty t^{s-1}\Theta_{\ba}(\tau;t)\dd t
 =\mu_{\ba}(-s)\Lambda_\tau(s).
\end{equation}
Hence the mixed Mellin transform has a meromorphic continuation via
\eqref{eq:negative-moment} and $\Lambda_\tau$, and its critical pole at $s=1$ has leading coefficient
$\mu_{\ba}(-1)$.
\end{theorem}

\begin{proof}
For real $s>0$, Tonelli's theorem and
$\Gamma(s)x^{-s}=\int_0^\infty y^{s-1}e^{-xy}\dd y$ first give
$\Gamma(s)\mu_{\ba}(-s)=\int_0^\infty y^{s-1}\Fa_{\ba}(y)\dd y$.
Using $y=x^2$ in this integral yields
\[
 \Gamma(s)\mu_{\ba}(-s)
 =2D\int_0^\infty
 x^{2s+d-1}\frac{e^{-Bx/2}}
 {\prod_{j=1}^d(1-e^{-|a_j|x})}\dd x,
\]
which is the standard Barnes-zeta integral.  In particular, the identity
shows that $\E(X_{\ba}^{-\sigma})<\infty$ for every real $\sigma>0$.
On each compact substrip of $\Re s>0$, the integrand $X_{\ba}^{-s}$ is
dominated by a sum of two finite negative moments, so
$s\mapsto\E(X_{\ba}^{-s})$ is holomorphic there.  The Barnes-zeta side is
holomorphic on the same half-plane; the identity theorem therefore proves
\eqref{eq:negative-moment} for all complex $s$ with $\Re s>0$.  For
\eqref{eq:Barnes-Epstein}, interchange the
absolutely convergent lattice sum and Mellin integral, then use
$y=\pi tQ_\tau(m,n)$ in each term.  The factor
$\int_0^\infty y^{s-1}\Fa_{\ba}(y)\dd y=\Gamma(s)\mu_{\ba}(-s)$ separates from the
Epstein sum.
\end{proof}

Combining the positive moments from Theorem~\ref{thm:master-probability}
with the negative moments just proved shows that
$\mu_{\ba}(z)=\E(X_{\ba}^z)$ is entire in $z$.  Indeed, on every compact
vertical strip $\alpha\le\Re z\le\beta$,
$|X_{\ba}^z|\le X_{\ba}^{\alpha}+X_{\ba}^{\beta}$, with the obvious
interpretation when an endpoint is zero, and both moments are finite.
Holomorphic dominated convergence therefore applies.  In particular,
$\mu_{\ba}$ is bounded on each closed vertical strip, a fact used in the
contour shifts below.

\begin{lemma}\label{lem:theta-endpoints-full}
For every $K>0$,
\begin{equation}\label{eq:theta-endpoints-full}
 \Theta_{\ba}(\tau;t)
 =\frac{\mu_{\ba}(-1)}{t}-1+O_K(t^K)
 \quad(t\downarrow0),
 \qquad
 \Theta_{\ba}(\tau;t)=O_K(t^{-K})
 \quad(t\to\infty).
\end{equation}
\end{lemma}

\begin{proof}
Let
\[
 \vartheta_\tau(u):=\sum_{(m,n)\in\Z^2}
 e^{-\pi uQ_\tau(m,n)}.
\]
Since $Q_\tau$ has determinant one, Poisson summation gives
$\vartheta_\tau(u)=u^{-1}\vartheta_\tau(u^{-1})$.  By the Laplace
representation and Tonelli,
\[
 \Theta_{\ba}(\tau;t)
 =\E\!\left[\vartheta_\tau(tX_{\ba})-1\right]
 =\frac{\mu_{\ba}(-1)}{t}-1
 +\frac1t\,\E\!\left[
 X_{\ba}^{-1}\bigl(\vartheta_\tau((tX_{\ba})^{-1})-1\bigr)
 \right].
\]
For every integer $L\ge1$, the theta transformation and the Gaussian
tail imply the global bound
\[
 \vartheta_\tau(y^{-1})-1\ll_{L,\tau} y^L,
 \qquad y>0.
\]
Indeed, for $0<y\le1$ this follows from the exponential decay of the theta
tail, while for $y\ge1$ the modular relation gives
$\vartheta_\tau(y^{-1})\ll_\tau y$.  Therefore the last term above is
\[
 \ll_{L,\tau} t^{L-1}\E[X_{\ba}^{L-1}],
\]
and choosing $L>K+1$ gives $O_K(t^K)$ by
Theorem~\ref{thm:master-probability}.

At infinity, for every integer $M\ge1$ we may retain $M$ factors in
\eqref{eq:Laplace-main} to obtain $\Fa_{\ba}(y)\ll_M(1+y)^{-M}$.
Choosing $M>\max\{K,1\}$ and using
$\sum_{\mathbf n\ne0}Q_\tau(\mathbf n)^{-M}<\infty$ gives
$\Theta_{\ba}(\tau;t)=O_K(t^{-K})$.
\end{proof}

\subsection{Classical Kronecker inputs and exact-order torsion}

The classical Kronecker formulas \cite{Kronecker1857} identify the geometric
quantities represented by the Mellin factors; the new input here is their
interaction with exact-order averaging.  For the determinant-one Epstein
factor of Theorem~\ref{thm:Barnes-Epstein-bridge}, the first Kronecker limit formula gives (see, e.g., \cite{Lang1987})
\begin{equation}\label{eq:Lambda-laurent-v122}
 \Lambda_\tau(s)=\frac1{s-1}+\kappa_\tau+O(s-1),
 \qquad
 \kappa_\tau=\gamma-\log(4\pi)
 -2\log\bigl(\sqrt v\,|\eta(\tau)|^2\bigr).
\end{equation}

For the second Kronecker formula fix the standard odd Jacobi theta
normalization $\vartheta_1'(0,\tau)=2\pi\eta(\tau)^3$.  For
$(\alpha,\beta)\notin\Z^2$ set
\begin{align}
 \chi_{\alpha,\beta}(m,n)&:=e^{2\pi i(m\alpha+n\beta)},\\
 Z_\tau(s;\alpha,\beta)
 &:=\sum_{(m,n)\ne(0,0)}
 \chi_{\alpha,\beta}(m,n)Q_\tau(m,n)^{-s},\\
 \Lambda_\tau(s;\alpha,\beta)
 &:=\pi^{-s}\Gamma(s)Z_\tau(s;\alpha,\beta),\\
 z_{\alpha,\beta}&:=\beta-\alpha\tau.
\end{align}
To avoid a normalization ambiguity, we distinguish the positive
multiplicative Arakelov Green function \cite{deJong2005} from its additive logarithmic version:
\begin{align}
 G_{\mathrm{Ar},\tau}(0,z)
 &:=\left|
 e^{-\pi(\operatorname{Im}z)^2/v}
 \frac{\vartheta_1(z,\tau)}{\eta(\tau)}\right|,
 \label{eq:Ar-Green-v122}\\
 \mathcal G_\tau(z)
 &:=-2\log G_{\mathrm{Ar},\tau}(0,z).
 \label{eq:additive-Green-v122}
\end{align}
These definitions fix all constants and signs used below; in particular,
$\mathcal G_\tau$ is the additive Green function associated with the displayed
multiplicative normalization, not a differently normalized Green kernel.
With precisely these conventions the classical second Kronecker limit formula reads
\begin{equation}\label{eq:second-KLF-v122}
 \Lambda_\tau(1;\alpha,\beta)
 =\mathcal G_\tau(z_{\alpha,\beta}).
\end{equation}
Equivalently,
\begin{equation}\label{eq:second-KLF-Z-v122}
 Z_\tau(1;\alpha,\beta)
 =2\pi^2\alpha^2v
 -2\pi\log\left|
 \frac{\vartheta_1(\beta-\alpha\tau,\tau)}{\eta(\tau)}\right|.
\end{equation}
See, for example, \cite{Lang1987}; the multiplicative normalization in
\eqref{eq:Ar-Green-v122} is compatible with the Arakelov normalization used
in \cite{deJong2005}.

Define the twisted mixed transform
\begin{equation}\label{eq:twisted-theta-v122}
 \Theta_{\ba}^{\alpha,\beta}(\tau;t)
 :=\sum_{(m,n)\ne(0,0)}
 \chi_{\alpha,\beta}(m,n)
 \Fa_{\ba}\!\left(\pi tQ_\tau(m,n)\right).
\end{equation}

The normalized remainder $\Ta_{m,\ba}$ from
Subsection~\ref{sec:remainders} has Mellin factor
\begin{equation}\label{eq:Bma-v12}
 B_{m,\ba}(s)
 :=\frac{(m+1)!\Gamma(1-s)}{\Gamma(m+2-s)}
 \frac{\mu_{\ba}(m+1-s)}{M_{m+1}(\ba)},
 \qquad
 A_{m,\ba}:=(m+1)\frac{M_m(\ba)}{M_{m+1}(\ba)}.
\end{equation}
Indeed, for $0<\Re s<1$, the Laplace representation
\eqref{eq:T-Laplace} and Tonelli--Fubini give
\begin{align*}
 \int_0^\infty y^{s-1}\Ta_{m,\ba}(y)\dd y
 &=\Gamma(s)\E[Z_{m,\ba}^{-s}]\\
 &=\Gamma(s)\E[U_m^{-s}]
   \E[(X_{\ba}^{[m+1]})^{-s}].
\end{align*}
Here
\[
 \E[U_m^{-s}]
 =\frac{(m+1)!\Gamma(1-s)}{\Gamma(m+2-s)},
 \qquad
 \E[(X_{\ba}^{[m+1]})^{-s}]
 =\frac{\mu_{\ba}(m+1-s)}{M_{m+1}(\ba)},
\]
so
\[
 \int_0^\infty y^{s-1}\Ta_{m,\ba}(y)\dd y
 =\Gamma(s)B_{m,\ba}(s).
\]
A direct two-dimensional lattice sum of $\Ta_{m,\ba}(\pi tQ_\tau)$ is not
absolutely convergent because $\Ta_{m,\ba}(y)\sim A_{m,\ba}/y$.
Accordingly we define the lattice normalized-remainder kernel by Mellin
regularization.  With
$\Lambda_\tau(s;0,0):=\Lambda_\tau(s)$, choose any $c\in(1,2)$ and set
\begin{equation}\label{eq:twisted-lattice-remainder}
 \mathcal R_{m,\ba}^{\alpha,\beta}(\tau;t)
 :=\frac{1}{2\pi i}\int_{c-i\infty}^{c+i\infty}
 B_{m,\ba}(s)\Lambda_\tau(s;\alpha,\beta)t^{-s}\dd s.
\end{equation}
The value is independent of $c\in(1,2)$: the integrand is holomorphic in
that open strip and has sufficient vertical decay for the standard contour
displacement argument.  We write $\mathcal R_{m,\ba}(\tau;t)$ for the
untwisted case.  This definition is the regularized theta lift naturally
associated with the normalized remainder;
after exact-order cancellation, the Mellin values at $s=1$ used below become
ordinary convergent logarithmic moments.

\begin{proposition}
\label{prop:mixed-second-KLF-v122}
For every nontrivial characteristic $(\alpha,\beta)\notin\Z^2$,
\begin{equation}\label{eq:twisted-endpoints-v122}
 \Theta_{\ba}^{\alpha,\beta}(\tau;t)
 =-1+O_K(t^K)\quad(t\downarrow0),
 \qquad
 \Theta_{\ba}^{\alpha,\beta}(\tau;t)=O_K(t^{-K})\quad(t\to\infty)
\end{equation}
for every $K>0$.  Hence the following is an ordinary convergent integral,
not merely a regularized value:
\begin{equation}\label{eq:mixed-second-KLF-v122}
 \int_0^\infty\Theta_{\ba}^{\alpha,\beta}(\tau;t)\dd t
 =\mu_{\ba}(-1)\mathcal G_\tau(\beta-\alpha\tau).
\end{equation}
\end{proposition}

\begin{proof}
Write the full twisted Gaussian theta sum as
\[
 \vartheta_{\tau}^{\alpha,\beta}(y)
 :=\sum_{(m,n)\in\Z^2}
 e^{2\pi i(m\alpha+n\beta)}e^{-\pi yQ_\tau(m,n)}.
\]
Poisson summation expresses it as $y^{-1}$ times a Gaussian sum over the
dual lattice shifted by $(\alpha,\beta)$.  Since
$(\alpha,\beta)\notin\Z^2$, the shifted dual lattice stays a positive
distance from the origin.  Consequently, for some $c=c(\tau,\alpha,\beta)>0$,
\[
 \vartheta_{\tau}^{\alpha,\beta}(y)
 \ll_{\tau,\alpha,\beta} y^{-1}e^{-c/y}
 \qquad(0<y\le1).
\]
In particular, for every integer $L\ge1$ one has the global bound
$\vartheta_{\tau}^{\alpha,\beta}(y)\ll_{L,\tau,\alpha,\beta}y^L$ for
$y>0$.  Using the Laplace representation
$\Fa_{\ba}(y)=\E(e^{-yX_{\ba}})$ and Tonelli therefore gives
\[
 \Theta_{\ba}^{\alpha,\beta}(\tau;t)+1
 =\E\!\left[\vartheta_{\tau}^{\alpha,\beta}(tX_{\ba})\right]
 \ll_{L,\tau,\alpha,\beta}
 t^L\E(X_{\ba}^L).
\]
All positive moments of $X_{\ba}$ are finite by
Theorem~\ref{thm:master-probability}; choosing $L>K$ proves the flat
expansion at the origin in \eqref{eq:twisted-endpoints-v122}.  Rapid decay
at infinity follows as in the untwisted case, since the character has unit
modulus and $\Fa_{\ba}(y)\ll_M(1+y)^{-M}$ for every $M$.

Termwise Mellin transformation is valid for $\Re s>1$ and gives
\[
 \int_0^\infty t^{s-1}\Theta_{\ba}^{\alpha,\beta}(\tau;t)\dd t
 =\mu_{\ba}(-s)\Lambda_\tau(s;\alpha,\beta).
\]
The endpoint estimates extend the defining integral holomorphically through
$s=1$.  Evaluating there and using \eqref{eq:second-KLF-v122} proves the
claim.
\end{proof}

\begin{proposition}
\label{prop:torsion-scale-v121}
For the twisted theta family above and the lattice normalized-remainder family in
\eqref{eq:twisted-lattice-remainder}, identify
$E_\tau[N]\cong(N^{-1}\Z/\Z)^2$ via
$(\alpha,\beta)\mapsto\beta-\alpha\tau$, and let the superscript $P$ denote
the corresponding additive characteristic.  Then, for every $N\ge1$ and $t>0$,
\begin{align}
 \sum_{P\in E_\tau[N]}\Theta_{\ba}^{P}(\tau;t)
 &=N^2\Theta_{\ba}(\tau;N^2t),\label{eq:torsion-scale-theta-v121}\\
 \sum_{P\in E_\tau[N]}\mathcal R_{m,\ba}^{P}(\tau;t)
 &=N^2\mathcal R_{m,\ba}(\tau;N^2t).\label{eq:torsion-scale-R-v121}
\end{align}
Consequently, if the superscript $\mathrm{prim}$ denotes exact order $N$,
then M\"obius inversion gives
\begin{align}
 \Theta_{\ba,N}^{\mathrm{prim}}(\tau;t)
 &=\sum_{d\mid N}\mu_{\rm M}(N/d)d^2
   \Theta_{\ba}(\tau;d^2t),\label{eq:prim-scale-theta-v121}\\
 \mathcal R_{m,\ba,N}^{\mathrm{prim}}(\tau;t)
 &=\sum_{d\mid N}\mu_{\rm M}(N/d)d^2
   \mathcal R_{m,\ba}(\tau;d^2t).\label{eq:prim-scale-R-v121}
\end{align}
Thus this is exactly the operator $\mathscr P_{N,2}$ from
\eqref{eq:abstract-P}, giving the analytic form of exact-order torsion
averaging for these two kernels.
\end{proposition}

\begin{proof}
Summing the additive characters over $(N^{-1}\Z/\Z)^2$ gives
\[
 \sum_{\alpha,\beta\in N^{-1}\Z/\Z}
 e^{2\pi i(m\alpha+n\beta)}
 =\begin{cases}
 N^2,&N\mid m\ \text{and}\ N\mid n,\\
 0,&\text{otherwise}.
 \end{cases}
\]
Reindexing the surviving lattice points yields
\eqref{eq:torsion-scale-theta-v121}.  The same character
orthogonality, now applied for $\Re s>1$ to the twisted completed Epstein
factors, gives
\[
 \sum_{P\in E_\tau[N]}\Lambda_\tau(s;P)
 =N^{2-2s}\Lambda_\tau(s).
\]
Multiplying by $B_{m,\ba}(s)$ and applying inverse Mellin transformation gives
\eqref{eq:torsion-scale-R-v121}.  Decomposing full $N$-torsion into
exact-order strata and applying M\"obius inversion gives
\eqref{eq:prim-scale-theta-v121}--\eqref{eq:prim-scale-R-v121}.
\end{proof}

The next identity is an exact-order M\"obius refinement of the classical torsion product formula of \cite{deJong2005}; we record it because it makes the ordinary von Mangoldt detector explicit.  Write
\[
 \Lambda(N):=\sum_{d\mid N}\mu_{\rm M}(N/d)\log d
\]
for the ordinary von Mangoldt function.

\begin{corollary}
\label{cor:primitive-green-v122}
For every $N\ge2$,
\begin{align}
 \sum_{0\ne P\in E_\tau[N]}\mathcal G_\tau(P)
 &=-2\log N,\label{eq:all-Green-v122}\\
 \sum_{P\in E_\tau[N]^{\mathrm{prim}}}\mathcal G_\tau(P)
 &=-2\Lambda(N).\label{eq:prim-Green-v122}
\end{align}
Equivalently,
\begin{equation}\label{eq:prim-theta-integral-v122}
 \int_0^\infty\Theta_{\ba,N}^{\mathrm{prim}}(\tau;t)\dd t
 =-2\mu_{\ba}(-1)\Lambda(N).
\end{equation}
In multiplicative Arakelov normalization this is
\begin{equation}\label{eq:prim-Ar-product-v122}
 \prod_{P\in E_\tau[N]^{\mathrm{prim}}}
 G_{\mathrm{Ar},\tau}(0,P)
 =e^{\Lambda(N)}
 =\begin{cases}
 p,&N=p^k,\\
 1,&N\text{ is not a prime power}.
 \end{cases}
\end{equation}
\end{corollary}

\begin{proof}
For a common cutoff $\varepsilon>0$, the scale law and subtraction of the
trivial characteristic give
\begin{align*}
 \sum_{0\ne P\in E_\tau[N]}
 \int_\varepsilon^\infty\Theta_{\ba}^{P}(\tau;t)\dd t
 &=\int_\varepsilon^\infty
 \bigl[N^2\Theta_{\ba}(\tau;N^2t)-\Theta_{\ba}(\tau;t)\bigr]\dd t\\
 &=-\int_\varepsilon^{N^2\varepsilon}
 \Theta_{\ba}(\tau;u)\dd u.
\end{align*}
Using
$\Theta_{\ba}(\tau;u)=\mu_{\ba}(-1)u^{-1}-1+O_K(u^K)$ and letting
$\varepsilon\downarrow0$ yields
$-2\mu_{\ba}(-1)\log N$.  Proposition~\ref{prop:mixed-second-KLF-v122}
then proves \eqref{eq:all-Green-v122}.  M\"obius inversion gives
\eqref{eq:prim-Green-v122} and \eqref{eq:prim-theta-integral-v122};
exponentiating \eqref{eq:prim-Green-v122} using
\eqref{eq:additive-Green-v122} gives \eqref{eq:prim-Ar-product-v122}.
The full-torsion product is also the classical identity
$\prod_{0\ne P\in E_\tau[N]}G_{\mathrm{Ar},\tau}(0,P)=N$ of
\cite{deJong2005}.
\end{proof}

\subsection{Arithmetic Mellin resonance on exact-order torsion}

Proposition~\ref{prop:torsion-scale-v121} identifies M\"obius inversion
from full torsion to exact order exactly with the abstract operator
$\mathscr P_{N,2}$ of \eqref{eq:abstract-P}.  Thus the elliptic critical
point $s=1$ is not merely analogous to the resonance point of
Section~\ref{sec:abstract-resonance}; it is a concrete realization of
Theorem~\ref{thm:AMR-local}.  In particular,
\begin{equation}\label{eq:torsion-abstract-identification}
 \widehat{\mathscr P_{N,2}H}(s)
 =J_{2-2s}(N)\widehat H(s),
\end{equation}
and each distinct prime divisor of $N$ removes one order of a critical
Mellin pole.

This interpretation also explains the primitive Green detector
\eqref{eq:prim-Green-v122}: the untwisted Epstein factor supplies the simple
base pole at $s=1$, while the second Kronecker limit formula identifies the
nontrivial-character values there with the additive Green function.  Hence the
critical layer is $\omega(N)=1$, and the first surviving coefficient is the
ordinary von Mangoldt function.  The normalized-remainder kernel below has a
double base pole and therefore moves the finite nonzero critical layer to
$\omega(N)=2$.

\paragraph{Theta and normalized-remainder layers.}
For the theta kernel of Theorem~\ref{thm:Barnes-Epstein-bridge}, the base Mellin factor
$\mu_{\ba}(-s)\Lambda_\tau(s)$ has a simple pole at $s=1$ with
$c_{-1}=\mu_{\ba}(-1)$.  For the regularized lattice remainder
\eqref{eq:twisted-lattice-remainder}, the factor $B_{m,\ba}$ in
\eqref{eq:Bma-v12} satisfies
\begin{equation}\label{eq:R-double-v12}
 B_{m,\ba}(s)\Lambda_\tau(s)
 =-\frac{A_{m,\ba}}{(s-1)^2}+O((s-1)^{-1}).
\end{equation}

\begin{proposition}
\label{prop:untwisted-R-endpoint-full}
The Mellin-regularized kernel \eqref{eq:twisted-lattice-remainder} has ordinary
Mellin transform
\begin{equation}\label{eq:untwisted-R-Mellin-full}
 \int_0^\infty t^{s-1}\mathcal R_{m,\ba}(\tau;t)\dd t
 =B_{m,\ba}(s)\Lambda_\tau(s),
 \qquad 1<\Re s<2.
\end{equation}
There is a constant $C_{m,\ba,\tau}$ such that, for every $K>0$,
\begin{equation}\label{eq:untwisted-R-endpoint-full}
 \mathcal R_{m,\ba}(\tau;t)
 =\frac{A_{m,\ba}}{t}\bigl(\log t+C_{m,\ba,\tau}\bigr)
 -1+O_K(t^K)
 \qquad(t\downarrow0).
\end{equation}
At the opposite endpoint,
\begin{equation}\label{eq:untwisted-R-infty-full}
 \mathcal R_{m,\ba}(\tau;t)=O(t^{-2})\qquad(t\to\infty),
\end{equation}
with faster decay when the pole at $s=2$ is absent.
\end{proposition}

\begin{proof}
The inverse Mellin integral in
\eqref{eq:twisted-lattice-remainder} converges on every vertical line
$1<c<2$: the completed Epstein factor has exponential gamma decay in
vertical strips, while $B_{m,\ba}$ has at most polynomial growth there.
We first use this inverse representation to obtain the endpoint bounds below.
Those bounds show that the resulting function has an ordinary Mellin transform
throughout $1<\Re s<2$; the standard Mellin inversion theorem then recovers
$B_{m,\ba}(s)\Lambda_\tau(s)$ on that strip, proving
\eqref{eq:untwisted-R-Mellin-full}.

Near $s=1$,
$B_{m,\ba}(s)=-A_{m,\ba}(s-1)^{-1}+O(1)$ and
$\Lambda_\tau(s)=(s-1)^{-1}+\kappa_\tau+O(s-1)$, so the product has a
double pole with leading coefficient $-A_{m,\ba}$.  Shifting the inverse
Mellin contour to the left across $s=1$ produces the
$t^{-1}\log t$ term and its accompanying $t^{-1}$ constant.  At $s=0$ one has $B_{m,\ba}(0)=1$, while the functional equation of the
completed determinant-one Epstein zeta function gives
$\operatorname*{Res}_{s=0}\Lambda_\tau(s)=-1$.  Hence the product has
residue $-1$ at $s=0$, giving the constant term $-1$.  Further left shifts
give $O_K(t^K)$.

For $t\to\infty$, shift the contour to the right.  The next possible pole is
$s=2$, coming from
$\Gamma(1-s)/\Gamma(m+2-s)
 =[(1-s)\cdots(m+1-s)]^{-1}$; hence the first possible term is $O(t^{-2})$.
For $m=0$ that pole is absent, and the contour can be shifted farther.
\end{proof}

\begin{proposition}
\label{prop:twisted-remainder-v122}
For the nontrivially twisted lattice normalized-remainder family defined in
\eqref{eq:twisted-lattice-remainder} and $1<\Re s<2$,
\begin{equation}\label{eq:twisted-R-Mellin-v122}
 \int_0^\infty t^{s-1}
 \mathcal R_{m,\ba}^{\alpha,\beta}(\tau;t)\dd t
 =B_{m,\ba}(s)\Lambda_\tau(s;\alpha,\beta).
\end{equation}
Since $\Lambda_\tau(s;\alpha,\beta)$ is holomorphic at $s=1$, the untwisted
double pole splits to a simple pole, and
\begin{equation}\label{eq:twisted-R-limit-v122}
 \lim_{t\downarrow0}t\,
 \mathcal R_{m,\ba}^{\alpha,\beta}(\tau;t)
 =-A_{m,\ba}\mathcal G_\tau(\beta-\alpha\tau).
\end{equation}
Thus the $t^{-1}\log t$ term disappears under every nontrivial additive
character.
\end{proposition}

\begin{proof}
Starting from the inverse Mellin definition
\eqref{eq:twisted-lattice-remainder}, contour shifting gives the endpoint
bounds stated below; these imply ordinary Mellin integrability for
$1<\Re s<2$.  The standard Mellin inversion theorem then recovers the
factor $B_{m,\ba}(s)\Lambda_\tau(s;\alpha,\beta)$ and proves
\eqref{eq:twisted-R-Mellin-v122}.  Near $s=1$,
$B_{m,\ba}(s)=-A_{m,\ba}(s-1)^{-1}+O(1)$, whereas
\eqref{eq:second-KLF-v122} gives
$\Lambda_\tau(1;\alpha,\beta)=\mathcal G_\tau(\beta-\alpha\tau)$.
Hence the residue at $s=1$ is
$-A_{m,\ba}\mathcal G_\tau(\beta-\alpha\tau)$.  Shifting the inverse
Mellin contour across this single pole gives
\eqref{eq:twisted-R-limit-v122}.  The next possible right-hand pole is at
$s=2$, so the twisted kernel is $O(t^{-2})$ at infinity.
\end{proof}

Hence Theorem~\ref{thm:AMR-local} gives two concrete exact-order layers.
For later reference define their meromorphic Mellin transforms by
\begin{align}
 \mathfrak M_{\Theta,N}(s)
 &:=J_{2-2s}(N)\mu_{\ba}(-s)\Lambda_\tau(s),
 \label{eq:Mtheta-def-full}\\
 \mathfrak M_{R,N}(s)
 &:=J_{2-2s}(N)B_{m,\ba}(s)\Lambda_\tau(s).
 \label{eq:MR-def-full}
\end{align}
Whenever the corresponding Mellin integral converges, these expressions agree
with the ordinary transforms of
$\Theta_{\ba,N}^{\mathrm{prim}}$ and
$\mathcal R_{m,\ba,N}^{\mathrm{prim}}$, respectively.

\begin{corollary}
\label{cor:two-kernels-v12}
If $\omega(N)=q\ge1$, then the exact-order theta Mellin transform satisfies
\begin{equation}\label{eq:theta-vM-v12}
 \mathfrak M_{\Theta,N}^{(q-1)}(1)
 =\frac{(-2)^q}{q}\,\mu_{\ba}(-1)\Lambda_q(N).
\end{equation}
If $\omega(N)=q\ge2$, then the exact-order normalized-remainder transform
satisfies
\begin{equation}\label{eq:R-vM-v12}
 \mathfrak M_{R,N}^{(q-2)}(1)
 =-\frac{(-2)^q}{q(q-1)}\,A_{m,\ba}\Lambda_q(N).
\end{equation}
In particular,
\begin{align}
 \omega(N)=1&:\quad
 \mathfrak M_{\Theta,N}(1)
 =-2\mu_{\ba}(-1)\Lambda(N),\label{eq:primepower-v12}\\
 N=p^ar^b,\ p\ne r&:\quad
 \mathfrak M_{R,N}(1)
 =-4A_{m,\ba}\log p\log r.\label{eq:twoprime-v12}
\end{align}
For $\omega(N)\ge3$ the remainder value at $s=1$ vanishes, but its first
surviving higher derivative is still given by \eqref{eq:R-vM-v12}.
\end{corollary}

\begin{proof}
For the theta layer, apply Theorem~\ref{thm:AMR-local} with $r=1$ and
$c_{-1}=\mu_{\ba}(-1)$.  This gives
\[
 (q-1)!\frac{(-2)^q}{q!}\Lambda_q(N)\mu_{\ba}(-1)
 =\frac{(-2)^q}{q}\mu_{\ba}(-1)\Lambda_q(N),
\]
which is \eqref{eq:theta-vM-v12}.  For the normalized-remainder layer,
\eqref{eq:R-double-v12} gives $r=2$ and $c_{-2}=-A_{m,\ba}$, so the same
theorem yields
\[
 (q-2)!\frac{(-2)^q}{q!}\Lambda_q(N)(-A_{m,\ba})
 =-\frac{(-2)^q}{q(q-1)}A_{m,\ba}\Lambda_q(N),
\]
which is \eqref{eq:R-vM-v12}.  The two displayed special cases follow from
$\Lambda_1(p^a)=\log p$ and
$\Lambda_2(p^ar^b)=2\log p\log r$ for $p\ne r$.
\end{proof}

\begin{corollary}
\label{cor:endpoint-cancellation-v122}
For every $N\ge2$ and every $K\ge1$,
\begin{align}
 \Theta_{\ba,N}^{\mathrm{prim}}(\tau;t)
 &=-J_2(N)+O_K(t^K),\label{eq:theta-endpoint-v122}\\
 \mathcal R_{m,\ba,N}^{\mathrm{prim}}(\tau;t)
 &=\frac{2A_{m,\ba}\Lambda(N)}{t}
 -J_2(N)+O_K(t^K).\label{eq:R-endpoint-v122}
\end{align}
Consequently, exact-order averaging removes the untwisted
$t^{-1}\log t$ term for every $N\ge2$, and the remaining $t^{-1}$ term is
present exactly on the prime-power layer.
\end{corollary}

\begin{proof}
Insert the endpoint expansions from Lemma~\ref{lem:theta-endpoints-full} and
Proposition~\ref{prop:untwisted-R-endpoint-full}, namely
\[
 \Theta_{\ba}(\tau;t)=\frac{\mu_{\ba}(-1)}t-1+O_K(t^K),
\]
and
\[
 \mathcal R_{m,\ba}(\tau;t)
 =\frac{A_{m,\ba}}t\bigl(\log t+C_{m,\ba,\tau}\bigr)
 -1+O_K(t^K)
\]
into \eqref{eq:prim-scale-theta-v121} and
\eqref{eq:prim-scale-R-v121}.  Since
$\sum_{d\mid N}\mu_{\rm M}(N/d)=0$ for $N\ge2$, the $t^{-1}$ theta term
and the common $t^{-1}(\log t+C_{m,\ba,\tau})$ remainder term cancel.
The surviving divisor sums are
$\sum_{d\mid N}\mu_{\rm M}(N/d)d^2=J_2(N)$ and
$\sum_{d\mid N}\mu_{\rm M}(N/d)\log d=\Lambda(N)$, which gives the formulas.
\end{proof}

\subsection{Analytic continuation versus genuine logarithmic Mellin moments}

The filtration above is a statement about meromorphic Mellin transforms.  A
finite value or derivative of the analytic continuation need not equal an
ordinary integral unless the corresponding logarithmic moment is integrable.
The following elementary lemma records the convention used in this section.

\begin{lemma}\label{lem:Mellin-log-v14}
Let $F:(0,\infty)\to\C$, $s_0\in\R$, and $m\ge0$.  Suppose there exists
$\eta>0$ such that
\begin{equation}\label{eq:log-integrability-v12}
 \int_0^\infty
 t^{s_0-1}|F(t)|\,(t^\eta+t^{-\eta})
 (1+|\log t|^m)\dd t<\infty.
\end{equation}
Then $\widehat F$ is holomorphic in the strip
$|\Re(s-s_0)|<\eta$ and
\begin{equation}\label{eq:log-moment-v12}
 \left.\frac{\dd^j}{\dd s^j}\widehat F(s)\right|_{s=s_0}
 =\int_0^\infty t^{s_0-1}F(t)(\log t)^j\dd t,
 \qquad 0\le j\le m.
\end{equation}
\end{lemma}

\begin{proof}
For $|\Re(s-s_0)|\le\eta/2$, one has
$t^{\Re(s-s_0)}\le t^{\eta/2}+t^{-\eta/2}$, and the same bound multiplied
by $|\log t|^j$ is dominated by the integrand in
\eqref{eq:log-integrability-v12}.  Holomorphic dominated convergence gives
the claim.
\end{proof}

\begin{corollary}
\label{cor:genuine-log-moments-v121}
For the Bernoulli--Barnes theta family, every exact-order kernel with $N\ge2$
has absolutely convergent logarithmic Mellin moments at $s=1$ of every
nonnegative degree.  The same is true for the normalized-remainder family
whenever $\omega(N)\ge2$.  Consequently the derivatives in
\eqref{eq:theta-vM-v12} and \eqref{eq:R-vM-v12} are genuine logarithmic
integrals in these ranges.
\end{corollary}

\begin{proof}
The exact torsion scale law together with the small-$t$ expansions retained
from the preceding theta/remainder analysis gives, for every $K\ge1$,
\[
 \Theta_{\ba,N}^{\mathrm{prim}}(\tau;t)
 =-J_2(N)+O(t^K),\qquad N\ge2,
\]
and
\[
 \mathcal R_{m,\ba,N}^{\mathrm{prim}}(\tau;t)
 =\frac{2A_{m,\ba}\Lambda(N)}{t}-J_2(N)+O(t^K).
\]
If $\omega(N)\ge2$, then $\Lambda(N)=0$, so the remainder kernel is also
bounded at the origin.  At infinity the exact-order theta kernel is rapidly
decreasing, while Proposition~\ref{prop:untwisted-R-endpoint-full} and the
finite divisor sum give
$\mathcal R_{m,\ba,N}^{\mathrm{prim}}(\tau;t)=O(t^{-2})$.
Thus the strengthened hypothesis of the lemma holds for every logarithmic
degree after choosing any sufficiently small $0<\eta<1$.
\end{proof}

\section{Dimension--arithmetic transition}
\label{sec:higher-dimensional-resonance}

We now prove the geometric transition announced in the introduction.  The
two-dimensional normalized-remainder transform in
Section~\ref{sec:arithmetic-mellin} has a double pole at $s=1$ because the
first inverse-power term of $\Ta_{m,\ba}$ collides with the Epstein pole.
In higher lattice dimension the Epstein pole moves to $s=\ell/2$, whereas
the scalar remainder factor $B_{m,\ba}(s)$ retains its poles at the integers
$1,\ldots,m+1$.  This produces a finite resonance set, precisely the even
dimensions $2,4,\ldots,2m+2$, at which the critical pole multiplicity
increases from one to two.  The mechanism may be summarized as
\[
 \begin{aligned}
 c_{m,k}(\ba)y^{-k}
 &\quad\Longleftrightarrow\quad
 B_{m,\ba}(s)\text{ has a simple pole at }s=k,\\
 \Lambda_Q(s)
 &\text{ has a simple pole at }s=\ell/2.
 \end{aligned}
\]
Hence, when $\ell=2k$,
\[
 \begin{aligned}
 B_{m,\ba}(s)\Lambda_Q(s)
 &\text{ has a double critical pole}\\
 &\Longrightarrow\quad
 \omega(N)=2\text{ is the finite nonzero critical layer}.
 \end{aligned}
\]
whereas in the nonresonant case the critical pole remains simple and the
finite nonzero layer is the prime-power layer $\omega(N)=1$.

Throughout this section let $\ell\ge1$, let $Q$ be a positive-definite real
$\ell\times\ell$ matrix with $\det Q=1$, and write
\[
 Q[\mathbf n]=\mathbf n^{\mathsf T}Q\mathbf n,
 \qquad \mathbf n\in\Z^\ell.
\]
The completed Epstein zeta function \cite{Epstein1903} is
\begin{equation}\label{eq:highdim-Epstein-v14}
 \Lambda_Q(s)
 :=\pi^{-s}\Gamma(s)
 \sum_{\mathbf n\in\Z^\ell\setminus\{0\}}Q[\mathbf n]^{-s}.
\end{equation}
It is meromorphic on $\C$ and has a simple pole at
$s_\ell:=\ell/2$ with residue $1$; see, for example, \cite{Terras1988}.

\subsection{Exact-order characters in higher dimension}

The abstract operator $\mathscr P_{N,\ell}$ has a direct lattice
interpretation.  For an absolutely convergent lattice kernel
\[
 H(t)=\sum_{\mathbf n\in\Z^\ell\setminus\{0\}}h(t,\mathbf n)
\]
and $\alpha\in\R^\ell/\Z^\ell$, put
\[
 H^\alpha(t)
 :=\sum_{\mathbf n\ne0}e^{2\pi i\langle\mathbf n,\alpha\rangle}
 h(t,\mathbf n).
\]
If $h(t,N\mathbf n)=h(N^2t,\mathbf n)$, finite Fourier orthogonality on
$(N^{-1}\Z^\ell)/\Z^\ell$ gives
\begin{equation}\label{eq:highdim-character-scale-v16}
 \sum_{\alpha\in(N^{-1}\Z^\ell)/\Z^\ell}H^\alpha(t)
 =N^\ell H(N^2t).
\end{equation}
M\"obius inversion from characters of level dividing $N$ to characters of
exact order $N$ is therefore exactly $\mathscr P_{N,\ell}$ from
\eqref{eq:abstract-P}.  All Mellin valuation statements in this section are
thus direct applications of Theorem~\ref{thm:AMR-local} at
$s_\ell=\ell/2$.

\subsection{Higher-dimensional Barnes--Epstein theta kernels}

Define
\begin{equation}\label{eq:highdim-theta-v14}
 \Theta_{\ba,Q}^{(\ell)}(t)
 :=\sum_{\mathbf n\in\Z^\ell\setminus\{0\}}
 \Fa_{\ba}(\pi tQ[\mathbf n]).
\end{equation}
For $N\ge1$ write
\begin{equation}\label{eq:highdim-theta-prim-v15}
 \Theta_{\ba,Q,N}^{(\ell),\mathrm{prim}}
 :=\mathscr P_{N,\ell}\Theta_{\ba,Q}^{(\ell)}.
\end{equation}
By the preceding orthogonality relation, this is the aggregate of the
additively twisted theta kernels over characters of exact order $N$.

\begin{theorem}
\label{thm:highdim-theta-v14}
Initially for $\Re s>\ell/2$,
\begin{equation}\label{eq:highdim-theta-mellin-v14}
 \int_0^\infty t^{s-1}\Theta_{\ba,Q}^{(\ell)}(t)\dd t
 =\mu_{\ba}(-s)\Lambda_Q(s).
\end{equation}
Hence the critical pole at $s=s_\ell$ is simple, with coefficient
\begin{equation}\label{eq:highdim-theta-residue-v14}
 c_{-1}^{(\Theta,\ell)}
 =\mu_{\ba}(-\ell/2)
 =2D\,\frac{\Gamma(\ell+d)}{\Gamma(\ell/2)}
 \zeta_d\!\left(\ell+d,\frac B2;\bb\right).
\end{equation}
Consequently, if $N=p^a$ is a prime power, the meromorphically continued critical value of
the exact-order transform is
\begin{equation}\label{eq:highdim-theta-prime-v14}
 \widehat{\Theta_{\ba,Q,N}^{(\ell),\mathrm{prim}}}(\ell/2)
 =-4D\,\frac{\Gamma(\ell+d)}{\Gamma(\ell/2)}
 \zeta_d\!\left(\ell+d,\frac B2;\bb\right)\log p.
\end{equation}
More generally, if $N\ge2$ and $q=\omega(N)$, its first surviving critical
derivative has degree $q-1$ and equals
\begin{equation}\label{eq:highdim-theta-first-jet-v14}
 (-2)^q(q-1)!\,
 \mu_{\ba}(-\ell/2)
 \prod_{p\mid N}\log p.
\end{equation}
\end{theorem}

\begin{proof}
Termwise Mellin transformation in the absolute-convergence half-plane gives
\eqref{eq:highdim-theta-mellin-v14}.  The residue statement follows from
$\operatorname*{Res}_{s=\ell/2}\Lambda_Q(s)=1$ and
\eqref{eq:negative-moment}.  Apply
Theorem~\ref{thm:AMR-local} with $r=1$ to obtain
\eqref{eq:highdim-theta-prime-v14} and
\eqref{eq:highdim-theta-first-jet-v14}.
\end{proof}

\subsection{Renormalized remainders and the resonance set}

Theorem~\ref{thm:normalized-remainder} gives
\begin{equation}\label{eq:T-inverse-coeff-v14}
 \Ta_{m,\ba}(y)
 =\sum_{k=1}^{m+1}c_{m,k}(\ba)y^{-k}+O_L(y^{-L}),
 \qquad L>m+1,
\end{equation}
with the coefficients $c_{m,k}(\ba)$ defined in
\eqref{eq:T-infty-cmk}.  In particular $c_{m,1}(\ba)=A_{m,\ba}$.

\begin{lemma}
\label{lem:subtracted-Mellin-v15}
For $0\le K\le m+1$ put
\[
 \Ta_{m,\ba}^{\langle K\rangle}(y)
 :=\Ta_{m,\ba}(y)-\sum_{k=1}^{K}c_{m,k}(\ba)y^{-k},
\]
with the sum empty for $K=0$.  If $K<m+1$, then its ordinary Mellin integral
converges absolutely for
\[
 K<\Re s<K+1,
\]
whereas for $K=m+1$ it converges absolutely for $\Re s>m+1$.
In either region,
\begin{equation}\label{eq:subtracted-Mellin-v15}
 \int_0^\infty y^{s-1}\Ta_{m,\ba}^{\langle K\rangle}(y)\dd y
 =\Gamma(s)B_{m,\ba}(s),
\end{equation}
where the right-hand side is understood by meromorphic continuation from
$0<\Re s<1$.
\end{lemma}

\begin{proof}
The expansion at the origin gives $\Ta_{m,\ba}(y)=1+O(y)$, while
\eqref{eq:T-inverse-coeff-v14} gives
$\Ta_{m,\ba}^{\langle K\rangle}(y)=O(y^{-K-1})$ at infinity when
$K<m+1$, and superalgebraic decay when $K=m+1$.  Hence the asserted strips
of absolute convergence follow directly.

To identify the Mellin transform, split the original integral at $1$.  The
standard continuation obtained from the first $K$ inverse-power terms is
\[
 \Gamma(s)B_{m,\ba}(s)
 =\int_0^1 y^{s-1}\Ta_{m,\ba}(y)\dd y
 +\int_1^\infty y^{s-1}\Ta_{m,\ba}^{\langle K\rangle}(y)\dd y
 +\sum_{k=1}^{K}\frac{c_{m,k}(\ba)}{k-s}.
\]
For $\Re s>K$ the last sum equals
$-\int_0^1 y^{s-1}\sum_{k=1}^{K}c_{m,k}(\ba)y^{-k}\dd y$.
Combining the two integrals over $(0,1)$ proves
\eqref{eq:subtracted-Mellin-v15}.  For $K=m+1$ the superalgebraic remainder
removes any upper restriction on $\Re s$.
\end{proof}

Put
\[
 K_{m,\ell}:=
 \min\!\left(m+1,\left\lfloor\frac\ell2\right\rfloor\right)
\]
and define the absolutely convergent renormalized lattice kernel
\begin{equation}\label{eq:highdim-R-v14}
 \mathcal R_{m,\ba,Q}^{(\ell)}(t)
 :=\sum_{\mathbf n\in\Z^\ell\setminus\{0\}}
 \left[
 \Ta_{m,\ba}(\pi tQ[\mathbf n])
 -\sum_{k=1}^{K_{m,\ell}}
 c_{m,k}(\ba)(\pi tQ[\mathbf n])^{-k}
 \right].
\end{equation}
For $K_{m,\ell}=0$ the inner sum is empty.

The following theorem identifies the resonance set completely and shows that
no further critical collisions occur.

\begin{theorem}
\label{thm:dimension-resonance-v14}
The Mellin transform of \eqref{eq:highdim-R-v14}, initially in a nonempty
strip of absolute convergence, has meromorphic continuation
\begin{equation}\label{eq:highdim-R-mellin-v14}
 \widehat{\mathcal R_{m,\ba,Q}^{(\ell)}}(s)
 =B_{m,\ba}(s)\Lambda_Q(s),
\end{equation}
where $B_{m,\ba}$ is defined in \eqref{eq:Bma-v12}.  At
$s=s_\ell=\ell/2$ its pole order is
\begin{equation}\label{eq:resonance-order-v14}
 r_{m,\ell}
 =1+\mathbf 1_{\{\ell\in\{2,4,\ldots,2m+2\}\}}.
\end{equation}
Thus the critical pole is double exactly for
$\ell=2k$ with $1\le k\le m+1$ and is simple otherwise.

If $\ell=2k\le2m+2$, the leading Laurent coefficient is
\begin{equation}\label{eq:resonance-leading-v14}
 c_{-2}^{(R,2k)}
 =-\frac{c_{m,k}(\ba)}{(k-1)!}.
\end{equation}
If $\ell$ is nonresonant, the leading coefficient is
\begin{equation}\label{eq:nonresonance-leading-v14}
 c_{-1}^{(R,\ell)}=B_{m,\ba}(\ell/2).
\end{equation}
\end{theorem}

\begin{proof}
Put $K=K_{m,\ell}$.  By Lemma~\ref{lem:subtracted-Mellin-v15}, the scalar
Mellin transform of the bracket in \eqref{eq:highdim-R-v14} is
$\Gamma(s)B_{m,\ba}(s)$.  If $K<m+1$, choose
\[
 \max\!\left(K,\frac\ell2\right)<\Re s<K+1;
\]
this strip is nonempty by the definition of $K$.  If $K=m+1$, choose instead
$\Re s>\max(m+1,\ell/2)$.  In either case the scalar Mellin integral and the
Epstein series are both absolutely convergent.  Fubini's theorem and the
change of variables $y=\pi tQ[\mathbf n]$ therefore give
\[
 \widehat{\mathcal R_{m,\ba,Q}^{(\ell)}}(s)
 =B_{m,\ba}(s)\,\pi^{-s}\Gamma(s)
  \sum_{\mathbf n\ne0}Q[\mathbf n]^{-s}
 =B_{m,\ba}(s)\Lambda_Q(s).
\]
Meromorphic continuation proves \eqref{eq:highdim-R-mellin-v14} globally.

Using the gamma quotient in \eqref{eq:Bma-v12},
\[
 B_{m,\ba}(s)
 =\frac{(m+1)!}{(1-s)(2-s)\cdots(m+1-s)}
 \frac{\mu_{\ba}(m+1-s)}{M_{m+1}(\ba)}.
\]
Since $\mu_{\ba}$ is entire and
$\mu_{\ba}(m+1-k)=M_{m+1-k}(\ba)>0$ for $1\le k\le m+1$, no cancellation
occurs at the zeros of the denominator.  Thus the only poles of
$B_{m,\ba}$ are the simple poles $1,\ldots,m+1$.  If $\ell/2$ is not one
of these integers, then the denominator in $B_{m,\ba}(\ell/2)$ is nonzero
and $\mu_{\ba}(m+1-\ell/2)=\E[X_{\ba}^{m+1-\ell/2}]>0$, so
$B_{m,\ba}(\ell/2)\ne0$.  The factor $\Lambda_Q$ has a simple pole,
with residue $1$, at $s=\ell/2$.  Hence the product has order two exactly
when $\ell/2\in\{1,\ldots,m+1\}$ and order one otherwise,
which proves \eqref{eq:resonance-order-v14}.

For completeness, at $s=k$ the denominator factorizes as
\[
 \prod_{j=1}^{m+1}(j-s)
 =(k-s)\!\prod_{\substack{1\le j\le m+1\\j\ne k}}(j-k)+O((s-k)^2),
\]
and
\[
 \prod_{j\ne k}(j-k)=(-1)^{k-1}(k-1)!(m+1-k)!.
\]
Therefore
\[
 \operatorname*{Res}_{s=k}B_{m,\ba}(s)
 =(-1)^k\frac{(m+1)!}{(k-1)!(m+1-k)!}
   \frac{M_{m+1-k}(\ba)}{M_{m+1}(\ba)}
 =-\frac{c_{m,k}(\ba)}{(k-1)!},
\]
where the last equality is \eqref{eq:T-infty-cmk}.  Multiplication by the
unit Epstein residue gives \eqref{eq:resonance-leading-v14}; in the
nonresonant case the value $B_{m,\ba}(\ell/2)$ itself is the leading
coefficient, proving \eqref{eq:nonresonance-leading-v14}.
\end{proof}

\subsection{Dimension--arithmetic transition}

The analytic pole collision of Theorem~\ref{thm:dimension-resonance-v14}
has a sharp arithmetic consequence.  Exact-order averaging contributes a
zero of order $\omega(N)$, so the change of the base pole multiplicity from
one to two forces the finite nonzero critical layer to move by one
prime-support level.  The next theorem makes this transition explicit.

Let
\begin{equation}\label{eq:highdim-R-prim-v14}
 \mathcal R_{m,\ba,Q,N}^{(\ell),\mathrm{prim}}
 :=\mathscr P_{N,\ell}\mathcal R_{m,\ba,Q}^{(\ell)}.
\end{equation}
For $\alpha\in\R^\ell/\Z^\ell$, define the twisted renormalized kernel
by inserting the factor
$e^{2\pi i\langle\mathbf n,\alpha\rangle}$ in the absolutely convergent
series \eqref{eq:highdim-R-v14}.  Finite Fourier orthogonality then identifies
\eqref{eq:highdim-R-prim-v14} with the aggregate over additive characters of
exact order $N$.

\begin{theorem}
\label{thm:dimension-arithmetic-v14}
Let $N\ge2$, $q=\omega(N)$, and let $r_{m,\ell}$ be given by
\eqref{eq:resonance-order-v14}.  Then
\begin{equation}\label{eq:phase-order-v14}
 \ord_{s=\ell/2}
 \widehat{\mathcal R_{m,\ba,Q,N}^{(\ell),\mathrm{prim}}}(s)
 =q-r_{m,\ell}.
\end{equation}
Consequently:
\begin{enumerate}[label=\textup{(\roman*)},leftmargin=2.5em]
\item If $\ell$ is nonresonant, the critical layer is the prime-power
 layer.  For $N=p^a$,
 \begin{equation}\label{eq:nonres-prime-v14}
  \widehat{\mathcal R_{m,\ba,Q,N}^{(\ell),\mathrm{prim}}}(\ell/2)
  =-2B_{m,\ba}(\ell/2)\log p.
 \end{equation}
 For $q\ge2$ the critical value vanishes, and the first surviving derivative
 has degree $q-1$.

\item If $\ell=2k\le2m+2$ is resonant, the critical layer is the
 two-prime layer.  For $N=p^ar^b$ with $p\ne r$,
 \begin{equation}\label{eq:res-two-prime-v14}
  \widehat{\mathcal R_{m,\ba,Q,N}^{(2k),\mathrm{prim}}}(k)
  =-\frac{4c_{m,k}(\ba)}{(k-1)!}\log p\log r.
 \end{equation}
 If $q\ge3$, the critical value vanishes and the first surviving derivative
 has degree $q-2$:
 \begin{equation}\label{eq:res-first-jet-v14}
  \left.\frac{\dd^{q-2}}{\dd s^{q-2}}
  \widehat{\mathcal R_{m,\ba,Q,N}^{(2k),\mathrm{prim}}}(s)
  \right|_{s=k}
  =(-1)^{q+1}2^q(q-2)!
   \frac{c_{m,k}(\ba)}{(k-1)!}
   \prod_{p\mid N}\log p.
 \end{equation}
\end{enumerate}
Thus the number of distinct primes required to reach a finite nonzero critical
value jumps from one to two exactly when the Epstein pole collides with one
of the finite inverse-power poles of the normalized remainder.
\end{theorem}

\begin{proof}
By \eqref{eq:highdim-R-mellin-v14}, exact-order averaging multiplies the base
transform by $J_{\ell-2s}(N)$.  Theorem~\ref{thm:AMR-local} therefore gives
\[
 \ord_{s=\ell/2}
 \widehat{\mathcal R_{m,\ba,Q,N}^{(\ell),\mathrm{prim}}}(s)
 =q-r_{m,\ell},
\]
which is \eqref{eq:phase-order-v14}.

Suppose first that $\ell$ is nonresonant.  Then $r_{m,\ell}=1$ and the
leading Laurent coefficient of the base transform is
$c_{-1}=B_{m,\ba}(\ell/2)$.  For $q=1$, formula
\eqref{eq:abstract-first-jet} has derivative order zero and gives
\[
 -2c_{-1}\Lambda_1(p^a)
 =-2B_{m,\ba}(\ell/2)\log p,
\]
which is \eqref{eq:nonres-prime-v14}.  If $q\ge2$, the valuation is positive,
so the critical value is zero and the first nonzero derivative has order
$q-1$.

Now let $\ell=2k\le2m+2$.  Then $r_{m,\ell}=2$ and, by
\eqref{eq:resonance-leading-v14},
$c_{-2}=-c_{m,k}(\ba)/(k-1)!$.  At the critical layer $q=2$,
\eqref{eq:abstract-first-jet} gives
\[
 \widehat{\mathcal R_{m,\ba,Q,N}^{(2k),\mathrm{prim}}}(k)
 =\frac{(-2)^2}{2!}\Lambda_2(N)c_{-2}.
\]
For $N=p^ar^b$ with $p\ne r$,
$\Lambda_2(N)=2\log p\log r$, and this is exactly
\eqref{eq:res-two-prime-v14}.  For $q\ge3$, the first nonzero derivative has
order $q-2$, and substituting
$\Lambda_q(N)=q!\prod_{p\mid N}\log p$ into
\eqref{eq:abstract-first-jet} yields
\eqref{eq:res-first-jet-v14}.  This also makes explicit why the finite
nonzero critical layer changes from $q=1$ to $q=2$ at the resonant
dimensions.
\end{proof}

\begin{remark}
Theorem~\ref{thm:dimension-arithmetic-v14} is formulated entirely in terms of
meromorphically continued Mellin transforms, so no endpoint integrability is
needed for its proof.  Whenever a concrete exact-order kernel satisfies a
strip-margin estimate through $s=\ell/2$, the corresponding meromorphically continued value
and its derivatives are ordinary Mellin integrals and logarithmic moments by
Lemma~\ref{lem:Mellin-log-v14}.  In particular, this is the interpretation
used in the two-dimensional critical cases proved in
Section~\ref{sec:arithmetic-mellin}.
\end{remark}

\section{Spectral hierarchy and normal crossing}
\label{sec:spectral-hierarchy}

The local transfer theorem becomes substantially more rigid when the Mellin
kernel comes from spectral geometry.  We first allow a general dilation order
$\rho>0$, then package all higher-depth kernels into one two-variable
holomorphic spectral family.  Its Mellin transform has two transverse simple
polar divisors; the complete sequence of higher pole orders is exactly the
parameter-jet shadow of this normal crossing.  Exact-order M\"obius inversion
then converts those jet orders into the full prime-support hierarchy.  The
ordinary heat trace realizes depth one, while the master family realizes all
depths two and higher from the same positive elliptic operator.  On a flat
torus the divisor transform is exactly the sum over torsion points of exact
order, so the analytic filtration acquires a direct geometric meaning.

\subsection{Exact-order dilation of arbitrary order}

Let $n\in\mathbb N$ and $\rho>0$.  For a Mellin kernel $H$ define
\begin{equation}\label{eq:spectral-general-P}
 (\mathscr P^{(\rho)}_{N,n}H)(t)
 :=\sum_{d\mid N}\mu_{\rm M}(N/d)d^n H(d^\rho t).
\end{equation}
The critical Mellin point for this scaling is
\[
 s_*:=\frac{n}{\rho}.
\]

\begin{proposition}
\label{prop:rho-transfer}
Assume that $\widehat H$ has a Laurent expansion
\[
 \widehat H(s_*+u)=\sum_{m=-r}^{\infty}c_mu^m,
 \qquad r\ge1,\quad c_{-r}\ne0.
\]
Then
\begin{equation}\label{eq:rho-factorization}
 \widehat{\mathscr P^{(\rho)}_{N,n}H}(s)
 =J_{n-\rho s}(N)\widehat H(s),
\end{equation}
and
\begin{equation}\label{eq:rho-Jordan-jet}
 J_{n-\rho(s_*+u)}(N)
 =\sum_{k\ge0}\frac{(-\rho)^k}{k!}\Lambda_k(N)u^k.
\end{equation}
If $q=\omega(N)$, then
\begin{equation}\label{eq:rho-valuation}
 \ord_{s=s_*}\widehat{\mathscr P^{(\rho)}_{N,n}H}=q-r.
\end{equation}
For $q\ge r$ the first surviving derivative is
\begin{equation}\label{eq:rho-first-jet}
 \left.\frac{\dd^{q-r}}{\dd s^{q-r}}
 \widehat{\mathscr P^{(\rho)}_{N,n}H}(s)\right|_{s=s_*}
 =(q-r)!(-\rho)^qc_{-r}\prod_{p\mid N}\log p.
\end{equation}
\end{proposition}

\begin{proof}
Termwise Mellin transformation in a common strip gives
\[
 \int_0^\infty t^{s-1}d^n H(d^\rho t)\dd t
 =d^{n-\rho s}\widehat H(s),
\]
which proves \eqref{eq:rho-factorization}.  Since $n-\rho s_* =0$,
Taylor expansion in $u$ gives \eqref{eq:rho-Jordan-jet}.  The identities
$\Lambda_k(N)=0$ for $k<q$ and
$\Lambda_q(N)=q!\prod_{p\mid N}\log p$ then show that the Jordan factor has
an exact zero of order $q$ with leading coefficient
$(-\rho)^q\prod_{p\mid N}\log p$.  Multiplication by the Laurent expansion of
$\widehat H$ proves \eqref{eq:rho-valuation} and
\eqref{eq:rho-first-jet}.
\end{proof}

\subsection{A two-variable spectral master family}

Let $M$ be a closed smooth manifold of dimension $n$, and let $P$ be a
positive self-adjoint classical elliptic operator of order $\rho>0$.  We
remove the kernel of $P$ and write its positive eigenvalues, repeated with
multiplicity, as $\{\lambda_\nu\}_{\nu\ge1}$.  Put
\[
 a:=\frac{n}{\rho},
 \qquad
 \zeta_P(s):=\operatorname{Tr}'P^{-s}
 =\sum_{\nu\ge1}\lambda_\nu^{-s}.
\]
By the classical theory of elliptic complex powers and Weyl asymptotics
\cite{Seeley1967,Hormander1968}, $\zeta_P$ extends meromorphically and its
rightmost pole is simple:
\begin{equation}\label{eq:spectral-zeta-leading}
 \zeta_P(a+x)=\frac{R_P}{x}+O(1),
 \qquad R_P>0.
\end{equation}

The depth-one member of the spectral hierarchy is supplied by the ordinary
heat trace
\[
 \Theta_P(t):=\operatorname{Tr}'e^{-tP}.
\]
For $\Re s>a$,
\begin{equation}\label{eq:spectral-heat-mellin}
 \widehat{\Theta_P}(s)=\Gamma(s)\zeta_P(s),
\end{equation}
so Proposition~\ref{prop:rho-transfer} gives
\begin{equation}\label{eq:spectral-heat-valuation}
 \ord_{s=a}\widehat{\mathscr P^{(\rho)}_{N,n}\Theta_P}(s)
 =\omega(N)-1.
\end{equation}
In particular, for $N=p^e$,
\begin{equation}\label{eq:spectral-heat-critical}
 \widehat{\mathscr P^{(\rho)}_{p^e,n}\Theta_P}(a)
 =-\rho\,\Gamma(a)R_P\log p.
\end{equation}

The higher-depth members arise simultaneously from a single holomorphic
parameter family.  For $u$ near $0$ and $y>0$, set
\begin{equation}\label{eq:spectral-master-scalar}
 \Psi_{a,u}(y)
 :=(1+y)^{-(a-u)}-y^{-(a-u)},
\end{equation}
and define
\begin{equation}\label{eq:spectral-master-trace}
 \mathfrak S_P(u;t)
 :=\sum_{\nu\ge1}\Psi_{a,u}(t\lambda_\nu)
 =\operatorname{Tr}'\!\left((I+tP)^{-(a-u)}-(tP)^{-(a-u)}\right).
\end{equation}
The cancellation in \eqref{eq:spectral-master-scalar} improves the
large-$y$ decay by one power.  The next lemma records the uniform trace-class
bounds and the Mellin-integrability estimates needed to justify all parameter
differentiations and interchanges below.

\begin{lemma}\label{lem:spectral-uniform-justification}
Fix $0<\eta<\min\{1,a\}$ and let
\[
 \Psi_{a,u}(y)=(1+y)^{-(a-u)}-y^{-(a-u)}.
\]
Let $U\Subset\{u\in\mathbb C:|u|<\eta\}$ and $J\ge0$.  Put
\[
 \beta_-:=\inf_{u\in U}\Re(a-u).
\]
Then, for $0\le j\le J$, uniformly in $u\in U$,
\begin{align}
 |\partial_u^j\Psi_{a,u}(y)|
 &\le C_{U,J}\,y^{-\Re(a-u)}(1+|\log y|^J),
 &&0<y\le1,\label{eq:spectral-uniform-small-y}\\
 |\partial_u^j\Psi_{a,u}(y)|
 &\le C_{U,J}\,y^{-\Re(a-u)-1}(1+(\log y)^J),
 &&y\ge1.\label{eq:spectral-uniform-large-y}
\end{align}
Consequently:
\begin{enumerate}[label=\textup{(\roman*)},leftmargin=2.5em]
\item for every $t>0$, the operator family
\[
 A_u(t):=(I+tP)^{-(a-u)}-(tP)^{-(a-u)}
\]
is holomorphic with values in the trace class on $|u|<\eta$, and
\[
 \partial_u^j\operatorname{Tr}'A_u(t)
 =\operatorname{Tr}'\partial_u^jA_u(t)
 =\sum_{\nu\ge1}\partial_u^j\Psi_{a,u}(t\lambda_\nu);
\]
\item on every compact subset of
\[
 \{(u,s):|u|<\eta,\ \max\{a,a-\Re u\}<\Re s<a-\Re u+1\},
\]
all sums and Mellin integrals obtained by taking at most $J$ derivatives in
$u$ are absolutely and locally uniformly convergent.  In particular,
Fubini--Tonelli and differentiation under both the trace and the Mellin
integral are valid there.
\end{enumerate}
\end{lemma}

\begin{proof}
Write $\beta=a-u$.  Since
\[
 \partial_u^j\Psi_{a,u}(y)
 =\frac{\log^j(1+y)}{(1+y)^\beta}
  -\frac{\log^j y}{y^\beta},
\]
the estimate \eqref{eq:spectral-uniform-small-y} follows uniformly on $U$.
For $y\ge1$ apply the fundamental theorem of calculus to
$h_{\beta,j}(x)=(\log x)^j x^{-\beta}$:
\[
 h_{\beta,j}(y+1)-h_{\beta,j}(y)
 =\int_y^{y+1}h'_{\beta,j}(x)\,\dd x.
\]
Because $u$ ranges over the compact set $U$, the derivative satisfies
\[
 |h'_{\beta,j}(x)|
 \le C_{U,J}x^{-\Re\beta-1}(1+(\log x)^J),
 \qquad x\ge1,
\]
and $x\asymp y$ on $[y,y+1]$.  This proves
\eqref{eq:spectral-uniform-large-y}.

Because $\eta<1$, one has $\beta_-+1>a$.  Weyl asymptotics therefore imply
\[
 \sum_{\nu\ge1}\lambda_\nu^{-\beta_- -1}
 (1+\log^J(2+\lambda_\nu))<\infty.
\]
For fixed $t>0$, \eqref{eq:spectral-uniform-large-y} and
$\Re(a-u)\ge\beta_-$ give a locally uniform trace-norm majorant for all but
finitely many eigenvalues, while the remaining terms are entire in $u$.
The Weierstrass theorem for Banach-valued holomorphic series proves (i).

For (ii), let $K$ be a compact subset of the displayed $(u,s)$-region and set
\[
 \delta_0:=\inf_K\bigl(\Re s-\Re(a-u)\bigr),\qquad
 \delta_1:=\inf_K\bigl(\Re(a-u)+1-\Re s\bigr),\qquad
 \delta_2:=\inf_K(\Re s-a).
\]
All three numbers are positive.  The bounds
\eqref{eq:spectral-uniform-small-y}--\eqref{eq:spectral-uniform-large-y}
therefore give, uniformly on $K$,
\[
 y^{\Re s-1}|\partial_u^j\Psi_{a,u}(y)|
 \ll_K y^{-1+\delta_0}(1+|\log y|^J),\qquad 0<y\le1,
\]
and
\[
 y^{\Re s-1}|\partial_u^j\Psi_{a,u}(y)|
 \ll_K y^{-1-\delta_1}(1+(\log y)^J),\qquad y\ge1.
\]
Thus the scalar Mellin integrals are absolutely and locally uniformly
convergent.  After the change of variables $y=t\lambda_\nu$, the remaining
eigenvalue factor is $\lambda_\nu^{-\Re s}$; the positive margin
$\delta_2$ and Weyl asymptotics give locally uniform convergence of the
eigenvalue sum.  Absolute Fubini and dominated convergence then justify all
claimed interchanges and differentiations.
\end{proof}

\begin{theorem}
\label{thm:spectral-master-family}
Fix $0<\eta<\min\{1,a\}$.  If $|u|<\eta$ and
\begin{equation}\label{eq:spectral-master-strip}
 \max\{a,a-\Re u\}<\Re s<a-\Re u+1,
\end{equation}
then
\begin{equation}\label{eq:spectral-master-Mellin}
 \widehat{\mathfrak S_P}(u,s)
 =\Gamma(s)\frac{\Gamma(a-u-s)}{\Gamma(a-u)}\zeta_P(s).
\end{equation}
The right-hand side gives a meromorphic continuation in $(u,s)$ near
$(0,a)$.  With $x=s-a$, its polar germ is a normal crossing of the two
simple divisors $x=0$ and $x+u=0$: more precisely, there are functions
$G,H,Z$ holomorphic near the origin, with $G(0,0)=1$, such that
\begin{align}
 \widehat{\mathfrak S_P}(u,a+x)
 &=-\frac{R_PG(u,x)}{x(x+u)}
   +\frac{R_PG(u,x)H(u,x)}{x}\notag\\
 &\hspace{2.5em}
   -\frac{G(u,x)Z(x)}{x+u}
   +G(u,x)Z(x)H(u,x).
 \label{eq:spectral-master-normal-crossing}
\end{align}
In particular,
\begin{equation}\label{eq:spectral-master-principal}
 \widehat{\mathfrak S_P}(u,a+x)
 =-\frac{R_P}{x(x+u)}
 +\text{terms with at most a simple pole along each divisor}.
\end{equation}

For $N\ge1$, set
\begin{equation}\label{eq:spectral-master-exact-order}
 \mathfrak S_{P,N}^{(\rho)}(u;t)
 :=\bigl(\mathscr P^{(\rho)}_{N,n}\mathfrak S_P(u;\cdot)\bigr)(t).
\end{equation}
If $q=\omega(N)$ and $W(N)=\prod_{p\mid N}\log p$, then, for every
$j\ge0$,
\begin{equation}\label{eq:spectral-master-jet-asymptotic}
 \left.\partial_u^j
 \widehat{\mathfrak S_{P,N}^{(\rho)}}(u,a+x)\right|_{u=0}
 =(-\rho)^qW(N)(-1)^{j+1}j!R_P\,
   x^{q-j-2}+O\!\left(x^{q-j-1}\right).
\end{equation}
Thus all support depths are encoded by the parameter jets of the single
meromorphic family \eqref{eq:spectral-master-Mellin}.
\end{theorem}

\begin{proof}
Put $\beta=a-u$.  For
$\Re\beta<\Re s<\Re\beta+1$,
\[
 (1+y)^{-\beta}-y^{-\beta}
 =-\beta\int_0^1(y+v)^{-\beta-1}\dd v.
\]
The beta integral and the change of variables $y=vw$ give
\begin{equation}\label{eq:spectral-master-scalar-Mellin}
 \int_0^\infty y^{s-1}\Psi_{a,u}(y)\dd y
 =\Gamma(s)\frac{\Gamma(a-u-s)}{\Gamma(a-u)}.
\end{equation}
Condition \eqref{eq:spectral-master-strip} simultaneously places the scalar
integral in this strip and gives $\sum_\nu\lambda_\nu^{-\Re s}<\infty$.
Lemma~\ref{lem:spectral-uniform-justification} supplies a locally uniform
absolute majorant on compact subsets of the strip.  After
$y=t\lambda_\nu$, Fubini--Tonelli therefore yields
\eqref{eq:spectral-master-Mellin}, and the same lemma gives holomorphic
dependence on $u$ in trace norm.

For the local polar structure, write
\[
 G(u,x):=\frac{\Gamma(a+x)}{\Gamma(a-u)},\qquad
 Z(x):=\zeta_P(a+x)-\frac{R_P}{x},
\]
and
\[
 H(u,x):=\Gamma(-u-x)+\frac1{u+x}.
\]
All three are holomorphic near $(0,0)$ and $G(0,0)=1$.  Substitution into
\eqref{eq:spectral-master-Mellin} gives
\eqref{eq:spectral-master-normal-crossing} exactly.

Finally, Proposition~\ref{prop:rho-transfer} gives
\[
 \widehat{\mathfrak S_{P,N}^{(\rho)}}(u,a+x)
 =J_{-\rho x}(N)\widehat{\mathfrak S_P}(u,a+x).
\]
Since
\[
 J_{-\rho x}(N)
 =(-\rho)^qW(N)x^q+O_N(x^{q+1}),
\]
the double-pole term in
\eqref{eq:spectral-master-normal-crossing} is the unique contribution that
can produce the power $x^{q-j-2}$ after $j$ derivatives in $u$.
Using
\[
 \left.\partial_u^j\frac1{x+u}\right|_{u=0}
 =(-1)^j j!x^{-j-1}
\]
proves \eqref{eq:spectral-master-jet-asymptotic}.
\end{proof}

The preceding theorem exhibits the hierarchy as the jet of a transverse
intersection of two simple polar divisors.  The next elementary local result
shows that this geometry is also the precise mechanism forcing the full
sequence of pole orders.

\begin{proposition}
\label{prop:normal-crossing-rigidity}
Let $F(u,x)$ be a meromorphic germ at $(0,0)$ of the form
\begin{equation}\label{eq:normal-crossing-general-germ}
 F(u,x)
 =\frac{A(u,x)}{x\,h(u,x)}
  +\frac{B(u,x)}{x}
  +\frac{C(u,x)}{h(u,x)}
  +D(u,x),
\end{equation}
where $A,B,C,D,h$ are holomorphic,
\[
 A(0,0)\ne0,\qquad h(0,0)=0,\qquad
 v:=\partial_xh(0,0)\ne0.
\]
Then $F(0,x)$ has a double pole at $x=0$.  Put
$c:=\partial_uh(0,0)$.  The following are equivalent:
\begin{enumerate}[label=\textup{(\roman*)},leftmargin=2.5em]
\item the divisors $x=0$ and $h(u,x)=0$ meet transversely at $(0,0)$;
\item $c\ne0$;
\item for every $j\ge0$, $\partial_u^jF(0,x)$ has a pole of exact order
      $j+2$ at $x=0$.
\end{enumerate}
When these conditions hold, the leading Laurent coefficient is
\begin{equation}\label{eq:normal-crossing-leading-jet}
 [x^{-j-2}]\,\partial_u^jF(0,x)
 =A(0,0)(-1)^j j!\frac{c^j}{v^{j+1}}.
\end{equation}
If $c=0$, the first parameter derivative $\partial_uF(0,x)$ has pole order
at most two, so the hierarchy already fails at $j=1$.
\end{proposition}

\begin{proof}
Because $v\ne0$, one has
$h(0,x)=vx+O(x^2)$, so the first term in
\eqref{eq:normal-crossing-general-germ} gives a genuine double pole and the
remaining terms have at most simple poles.  The two divisors are transverse
exactly when $dx$ and $dh$ are linearly independent at the origin, which is
equivalent to $c\ne0$.

Assume $c\ne0$.  In $\partial_u^j(h^{-1})$, the unique term with denominator
$h^{j+1}$ is
\[
 (-1)^j j!(\partial_uh)^j h^{-j-1};
\]
all other terms have smaller powers of $h$ in the denominator.  Hence the
first summand in \eqref{eq:normal-crossing-general-germ} has leading term
\[
 A(0,0)(-1)^j j!\frac{c^j}{v^{j+1}}x^{-j-2}.
\]
The $B/x$ term has pole order at most one, the $C/h$ term has pole order at
most $j+1$, and $D$ is holomorphic, so no cancellation of the displayed
coefficient is possible.  This proves (iii) and
\eqref{eq:normal-crossing-leading-jet}.

Conversely, if $c=0$, then $\partial_uh(0,x)=O(x)$.  Differentiating
\eqref{eq:normal-crossing-general-germ} once shows that the potentially most
singular term
$-A(0,x)\partial_uh(0,x)/(x h(0,x)^2)$ is only $O(x^{-2})$.
All other terms also have pole order at most two.  Thus (iii) fails for
$j=1$.
\end{proof}

For $j\in\mathbb Z_{\ge0}$ define
\begin{equation}\label{eq:Bessel-scalar-kernel}
 \Phi_{a,j}(y)
 :=\left.\partial_u^j\Psi_{a,u}(y)\right|_{u=0}
 =\frac{\log^j(1+y)}{(1+y)^a}
   -\frac{\log^j y}{y^a},
 \qquad y>0,
\end{equation}
and
\begin{equation}\label{eq:spectral-Bessel-trace}
 \mathcal S_{P,j}(t)
 :=\left.\partial_u^j\mathfrak S_P(u;t)\right|_{u=0}
 =\sum_{\nu\ge1}\Phi_{a,j}(t\lambda_\nu).
\end{equation}
Equivalently,
\[
 \mathcal S_{P,j}(t)
 =\left.(-\partial_\beta)^j
 \operatorname{Tr}'\!\left((I+tP)^{-\beta}-(tP)^{-\beta}\right)
 \right|_{\beta=a}.
\]
The endpoint estimates
\begin{align}
 \Phi_{a,j}(y)
 &=O_{a,j}\!\left(y^{-a}(1+|\log y|^j)\right)
 &&(y\downarrow0),\label{eq:Bessel-scalar-origin}\\
 \Phi_{a,j}(y)
 &=O_{a,j}\!\left(y^{-a-1}(1+|\log y|^j)\right)
 &&(y\to\infty)\label{eq:Bessel-scalar-decay}
\end{align}
show directly that the series in \eqref{eq:spectral-Bessel-trace} converges
absolutely for every $t>0$.

\begin{corollary}
\label{thm:spectral-support-hierarchy}
For every $j\ge0$, initially in the strip $a<\Re s<a+1$,
\begin{equation}\label{eq:spectral-Bessel-factorization}
 \widehat{\mathcal S_{P,j}}(s)
 =\left.\partial_u^j\left[
 \Gamma(s)\frac{\Gamma(a-u-s)}{\Gamma(a-u)}\zeta_P(s)
 \right]\right|_{u=0},
\end{equation}
and the right-hand side gives its meromorphic continuation.  At $s=a$ it
has a pole of exact order $j+2$ with leading Laurent coefficient
\begin{equation}\label{eq:spectral-leading-coeff}
 c_{-(j+2)}^{(P,j)}=(-1)^{j+1}j!R_P.
\end{equation}
Consequently, if $q=\omega(N)$,
\begin{equation}\label{eq:spectral-support-valuation}
 \ord_{s=a}
 \widehat{\mathscr P^{(\rho)}_{N,n}\mathcal S_{P,j}}(s)
 =q-(j+2).
\end{equation}
On the first finite nonzero layer $q=j+2$,
\begin{equation}\label{eq:spectral-critical-value}
 \widehat{\mathscr P^{(\rho)}_{N,n}\mathcal S_{P,j}}(a)
 =-\rho^{j+2}j!R_P\prod_{p\mid N}\log p.
\end{equation}
More generally, for $q\ge j+2$ the first nonzero derivative is
\begin{equation}\label{eq:spectral-first-derivative}
 \left.\frac{\dd^{q-j-2}}{\dd s^{q-j-2}}
 \widehat{\mathscr P^{(\rho)}_{N,n}\mathcal S_{P,j}}(s)
 \right|_{s=a}
 =(q-j-2)!(-\rho)^q(-1)^{j+1}j!R_P
 \prod_{p\mid N}\log p.
\end{equation}
\end{corollary}

\begin{proof}
Differentiate Theorem~\ref{thm:spectral-master-family} $j$ times in $u$ at
$u=0$.  Lemma~\ref{lem:spectral-uniform-justification}, applied on a compact
parameter disc and on compact substrips of $a<\Re s<a+1$, gives absolute
local uniform convergence of the differentiated trace series and Mellin
integrals.  It therefore justifies differentiation termwise under the trace,
summation, and Mellin integration.  Thus
\eqref{eq:spectral-master-normal-crossing} and
Proposition~\ref{prop:normal-crossing-rigidity} give the exact pole order and
\eqref{eq:spectral-leading-coeff}.  Formula
\eqref{eq:spectral-master-jet-asymptotic} then gives
\eqref{eq:spectral-support-valuation}--\eqref{eq:spectral-first-derivative}.
\end{proof}

\begin{corollary}\label{cor:complete-spectral-hierarchy}
For a fixed positive self-adjoint elliptic operator $P$ as above, every
positive critical prime-support depth is realized.  Depth one is realized by
the heat trace $\Theta_P$, while for each $r\ge2$ depth $r$ is realized by
$\mathcal S_{P,r-2}$.  Thus the same spectral datum $P$ produces the full
hierarchy
\[
 \omega(N)=1,2,3,\ldots.
\]
Moreover, the entire hierarchy is the parameter-jet shadow of the single
normal crossing $x(x+u)=0$ in
\eqref{eq:spectral-master-normal-crossing}.
\end{corollary}

\section{Flat tori and exact-order torsion}
\label{sec:flat-tori}

Let $Q$ be a positive definite $n\times n$ matrix with $\det Q=1$, and write
$\mathbb T^n=\mathbb R^n/\mathbb Z^n$.  For
$e_{\mathbf k}(x)=e^{2\pi i\langle\mathbf k,x\rangle}$ define the scaled
flat-torus Laplacian $P_Q$ by
\[
 P_Qe_{\mathbf k}=\pi Q[\mathbf k]e_{\mathbf k},
 \qquad \mathbf k\ne0.
\]
Thus $P_Q$ has order $2$, $a=n/2$, and
\begin{equation}\label{eq:flat-zeta-completed}
 \zeta_{P_Q}(s)
 =\pi^{-s}\sum_{\mathbf k\ne0}Q[\mathbf k]^{-s}
 =\frac{\Lambda_Q(s)}{\Gamma(s)}.
\end{equation}
With the determinant-one normalization used in
Section~\ref{sec:higher-dimensional-resonance},
$\operatorname*{Res}_{s=n/2}\Lambda_Q(s)=1$; hence
\begin{equation}\label{eq:flat-spectral-residue}
 R_{P_Q}=\frac1{\Gamma(n/2)}.
\end{equation}

Define the absolutely convergent Fourier kernel
\begin{equation}\label{eq:flat-Bessel-kernel}
 \mathcal K_{Q,j}(x;t)
 :=\sum_{\mathbf k\in\mathbb Z^n\setminus\{0\}}
 e_{\mathbf k}(x)
 \Phi_{n/2,j}\!\left(\pi tQ[\mathbf k]\right).
\end{equation}
For $N\ge1$, let $\mathbb T^n[N]^{\mathrm{ex}}$ be the points of exact
additive order $N$, and set
\begin{equation}\label{eq:flat-exact-torsion-aggregate}
 \mathcal K_{Q,j,N}^{\mathrm{ex}}(t)
 :=\sum_{x\in\mathbb T^n[N]^{\mathrm{ex}}}\mathcal K_{Q,j}(x;t).
\end{equation}

\begin{theorem}
\label{thm:flat-torus-hierarchy}
For every $N\ge1$,
\begin{equation}\label{eq:flat-torsion-Mobius}
 \mathcal K_{Q,j,N}^{\mathrm{ex}}(t)
 =\sum_{d\mid N}\mu_{\rm M}(N/d)d^n
  \mathcal K_{Q,j}(0;d^2t)
 =\bigl(\mathscr P_{N,n}\mathcal K_{Q,j}(0;\cdot)\bigr)(t).
\end{equation}
If $q=\omega(N)$, then
\begin{equation}\label{eq:flat-torsion-valuation}
 \ord_{s=n/2}\widehat{\mathcal K_{Q,j,N}^{\mathrm{ex}}}(s)
 =q-(j+2).
\end{equation}
The finite nonzero critical layer is therefore exactly $q=j+2$, and on that
layer
\begin{equation}\label{eq:flat-torsion-critical}
 \widehat{\mathcal K_{Q,j,N}^{\mathrm{ex}}}(n/2)
 =-\frac{2^{j+2}j!}{\Gamma(n/2)}
 \prod_{p\mid N}\log p.
\end{equation}
Thus the logarithmic hierarchy realizes successively the two-prime,
three-prime, four-prime, and all higher prime-support layers.  The depth-one
member is the ordinary heat kernel
\[
 \mathcal H_Q(x;t):=\sum_{\mathbf k\ne0}
 e_{\mathbf k}(x)e^{-\pi tQ[\mathbf k]}.
\]
Its exact-order torsion aggregate satisfies the same M\"obius identity as
\eqref{eq:flat-torsion-Mobius}; for $N=p^e$ its critical Mellin value is
$-2\log p$.  Hence flat-torus spectral geometry realizes every positive
critical prime-support depth by exact-order torsion summation.
\end{theorem}

\begin{proof}
For the full $d$-torsion subgroup, finite Fourier orthogonality gives
\[
 \sum_{x\in\mathbb T^n[d]}e_{\mathbf k}(x)
 =d^n\,\mathbf1_{d\mid\mathbf k},
\]
where $d\mid\mathbf k$ means divisibility of every coordinate.  Since
$Q[d\mathbf m]=d^2Q[\mathbf m]$, absolute convergence of
\eqref{eq:flat-Bessel-kernel} permits termwise summation and yields
\[
 \sum_{x\in\mathbb T^n[d]}\mathcal K_{Q,j}(x;t)
 =d^n\mathcal K_{Q,j}(0;d^2t).
\]
Every $d$-torsion point has exact order dividing $d$.  M\"obius inversion on
the divisor lattice therefore proves \eqref{eq:flat-torsion-Mobius}.
The function $\mathcal K_{Q,j}(0;t)$ is precisely
$\mathcal S_{P_Q,j}(t)$.  Apply Theorem~\ref{thm:spectral-support-hierarchy}
with $\rho=2$ and use \eqref{eq:flat-spectral-residue}; this gives
\eqref{eq:flat-torsion-valuation} and \eqref{eq:flat-torsion-critical}.
For $\mathcal H_Q$ the same Fourier orthogonality gives the corresponding
exact-order torsion identity, while
\[
 \widehat{\mathcal H_Q(0;\cdot)}(s)=\Lambda_Q(s).
\]
Since $\operatorname*{Res}_{s=n/2}\Lambda_Q(s)=1$, the depth-one formula is
Proposition~\ref{prop:rho-transfer} with $\rho=2$, $r=1$, and $c_{-1}=1$,
which gives $-2\log p$ on every prime-power layer $N=p^e$.
\end{proof}

\section{Conclusion}

Exact-order M\"obius inversion turns a Mellin pole of order $r$ into the
critical valuation $\omega(N)-r$, and its first surviving Jordan--Mellin jet
selects the weight $W(N)=\prod_{p\mid N}\log p$.  This same weight governs the
moving boundary law.  With several Galois splitting marks, the resulting
entanglement polynomial simultaneously controls mixed cumulants, local moving
profiles, and large deviations; its coefficients recover reciprocal compositum
degrees, while normalized pairwise covariance recovers field-intersection
degrees.

The spectral realization shows that the mechanism is structural rather than
kernel-specific.  A single two-variable spectral difference family attached to any
positive self-adjoint elliptic operator has a transverse normal crossing of
Mellin poles whose parameter jets generate all higher pole orders; exact-order
dilation converts them into every positive prime-support depth.  On flat tori
this is genuine exact-order torsion summation.  The Bernoulli--Barnes/
Barnes--Epstein construction provides a complementary explicit model in which
a geometric pole collision raises the critical multiplicity and shifts the
finite nonzero layer from one-prime to two-prime support.  Together these
results identify pole multiplicity, exact-order inversion, and marked prime
support as parts of one transfer mechanism across arithmetic, probabilistic,
spectral, and geometric settings.

\subsection*{Open problems}
The unified theorem suggests several directions in which the same mechanism
should persist but is not yet covered by the present paper.

\begin{problem}[Multiparameter transfer beyond a normal crossing]
Develop an exact-order transfer theorem for genuinely multiparameter Mellin
families in which several polar divisors meet and the mixed polar terms do not
factor into independent one-variable contributions.  The main issue is to
identify the resulting multi-valuation and the first surviving mixed jet in
terms of the arithmetic support semigroup.
\end{problem}

\begin{problem}[Higher Galois reconstruction]
Determine which compositum, multiple-intersection, or related lattice
invariants of a family of finite Galois extensions are determined by the full
collection of mixed prime-support cumulants, especially when the Frobenius
marking lattice is not of full rank.
\end{problem}

\begin{problem}[Beyond flat exact-order torsion]
Find geometric settings beyond flat tori in which the full spectral
normal-crossing hierarchy admits a direct exact-order geometric
interpretation, and determine whether higher-codimension polar intersections
produce new arithmetic support filtrations.
\end{problem}

\section*{Declarations}
\paragraph{\bf Authors' Contributions}
Both authors contributed to the mathematical development and preparation of the manuscript and approved the final version.

\paragraph{\bf Funding}
No funding was received for this research.

\paragraph{\bf Availability of data and materials}
Data sharing is not applicable to this article because no new data were created
or analyzed in this study.

\paragraph{\bf Competing interests}
The authors declare that they have no competing interests.

\end{document}